\documentclass[11pt]{article}
\usepackage{amsfonts,latexsym,amsmath,amscd,geometry}
\usepackage{amssymb}
\usepackage{latexsym}

\newcommand \nc{\newcommand}
\newtheorem{theorem}{Theorem}[section]
\newtheorem{lemma}[theorem]{Lemma}

\newtheorem{remark}[theorem]{Remark}

\nc{\ba}{\begin{array}}\nc{\ea}{\end{array}}
\nc{\be}{\begin{eqnarray}}\nc{\ee}{\end{eqnarray}}
\nc{\beq}{\begin{equation}}\nc{\eeq}{\end{equation}}
\nc{\bex}{\begin{eqnarray*}}\nc{\eex}{\end{eqnarray*}}
\nc{\btm}{\begin{theorem}} \nc{\etm}{\end{theorem}}
\nc{\blm}{\begin{lemma}} \nc{\elm}{\end{lemma}}
\nc{\R}{\mathbb{R}}  \nc{\ld}{\lambda}
\nc{\va}{\varphi}
\nc{\ve}{\varepsilon}
\def\x{\mathbf{x}}\def\e{\mathbf{e}}

\def\pa{\partial}
\def\pf{\noindent{\bf Proof.\quad}}\def\endpf{\hfill$\Box$}

\def\D{\widetilde{D}}\def\del{\widetilde{\Delta}}

\newcommand \qed {\hfill $\Box$}

\begin{document}
\title{Finite time singularity of the nematic liquid crystal flow in dimension three}
\author{Tao Huang\footnote{Department of Mathematics, The Pennsylvania State University, University Park, PA 16802, USA}\quad{\ Fanghua Lin\footnote{Courant Institute of Mathematical Sciences, New York University, NY 10012, USA}} \quad Chun Liu$^{*}$\quad Changyou Wang\footnote{Department of Mathematics, Purdue University, 150 N. University Street, West Lafayette, IN 47907, USA} }
\maketitle

\begin{abstract} In this paper, we consider the initial and boundary value problem of
a simplified nematic liquid crystal flow in dimension three and construct two examples of finite time
singularity. The first example is constructed within the class of axisymmetric solutions, while the second
example is constructed for any generic initial data $(u_0,d_0)$ that has sufficiently small energy, and
$d_0$ has a nontrivial topology.
\end{abstract}

\section {Introduction}
\setcounter{equation}{0}
\setcounter{theorem}{0}
Let $\Omega\subset\R^n\ (n=2, 3)$ be a bounded, smooth domain, and $0<T\le +\infty$.
In this paper, we will consider a simplified version of the hydrodynamic flow of nematic liquid crystals on $\Omega\times (0,T)$ given by
\begin{equation}\label{liquidcrystal}
\begin{cases}
u_t+u\cdot\nabla u-\mu\Delta u+\nabla P=-\lambda\nabla\cdot\big(\nabla d\odot\nabla d-\frac12{|\nabla d|^2}\mathbb{I}_n\big),\\
\nabla \cdot u=0, \\
d_t+u\cdot\nabla d=\gamma\big(\Delta d+|\nabla d|^2d\big),
\end{cases}
\end{equation}
where $u(\x,t):\Omega\times(0,T)\rightarrow\R^n$ is the velocity field of the underlying incompressible fluid,
$d(\x,t):\Omega\times(0,T)\rightarrow \mathbb S^2:=\big\{v\in\R^3:  |v|=1\big\}$ represents the (averaged) orientation field of nematic liquid crystal molecules,
$P(\x,t):\Omega\times(0,T)\rightarrow\R$ is the pressure function, $\x\in\Omega$,
$\nabla\cdot$ denotes the divergence operator on $\mathbb R^n$, $\nabla
d\odot\nabla d =\big(\langle\frac{\partial d}{\partial
\x_i},\frac{\partial d}{\partial \x_j}\rangle\big)_{1\leq i,j\leq n}\in \R^{n\times n}$
represents the stress tensor induced by the orientation field $d$,
and $\mathbb{I}_n=\big(\delta_{ij}\big)_{1\le i, j\le n}\in\R^{n\times n}$ is the identity matrix of order $n$. 
The parameters $\mu$, $\lambda$ and $\gamma$ are positive constants representing the fluid viscosity,
the competition between kinetic energy and potential energy, and the macroscopic elastic relaxation time for the molecular orientation field respectively.

The system \eqref{liquidcrystal}, first proposed by Lin \cite{lin},  is a simplified version of the general Ericksen-Leslie system modeling the hydrodynamic flow of nematic
liquid crystal materials proposed by Ericksen \cite{ericksen} and Leslie
\cite{leslie} during the period between 1958 and 1968. The system (\ref{liquidcrystal}) is a macroscopic continuum description of the time evolution of the material under the influence of both
the fluid velocity field and the macroscopic description of the microscopic
orientation configurations of rod-like liquid crystals.
The interested readers can refer to
\cite{ericksen}, \cite{leslie}, \cite{lin}, and Lin-Liu \cite{lin-liu} for
more details.

The system \eqref{liquidcrystal} is a strongly coupling system between the incompressible Naiver-Stokes equation and the
heat flow of harmonic maps into $\mathbb S^2$, and relates to several important equations:
\begin{itemize}
\item[1.] When $\lambda\equiv 0$ or $d\equiv e_0\in\mathbb S^2$, the system \eqref{liquidcrystal}$_{1,2}$
reduces to the incompressible Navier-Stokes equation (or NSE) that has been extensively studied for decays (see Lions \cite{lions}, Temam \cite{temam}). Although the existence of global weak solutions to the initial value problem to NSE
has been established by Leray \cite{leray} in 1930's and Hopf \cite{hopf} in 1950's,
it is a long outstanding open problem whether NSE  admits a global smooth
solution for a smooth initial data in dimension $n=3$. It is also an open problem
whether Leray-Hopf weak solutions to NSE are unique in dimension $n=3$.

\item[2.] When $u\equiv 0$ and $\gamma=1$, the equation \eqref{liquidcrystal}$_3$ reduces to the heat flow of harmonic maps into $\mathbb S^2$. For dimension $n=2$, the existence of a unique global weak solution, which has at most
finitely many singular points, has been proved by Struwe \cite{struwe} and Chang \cite{chang}.
In higher dimensions, the existence of a global, partially regular weak solution
has also been obtained by  Chen-Struwe \cite{chen-struwe} and Chen-Lin \cite{chen-lin}.
Examples of finite time singularities have been constructed by Coron-Ghidaglia \cite{cg}
and Chen-Ding \cite{chen-ding} for $n\ge 3$. In an important work \cite{chang-ding-ye}, Chang-Ding-Ye constructed
examples of finite time singularities when $n=2$ (see Grotowaski \cite{grotowaski, grotowaski2} for some generalizations to $n\ge 3$) by studying the equation
\begin{equation}\label{shm}
\va_t=\va_{rr}+\frac{\va_r}{r}-\frac{\sin(2\va)}{2r^2}, \ (r,t)\in (0,1)\times (0,+\infty).
\end{equation}
The interested readers can refer to Lin-Wang \cite{lin-wangbk} and references therein for more details.

\item[3.] Another important case we want to mention is $\gamma\equiv 0$. The system \eqref{liquidcrystal}
for $(u,d)$  is closely related to the MHD system for $(u,\psi)$
provided we identify $\psi=\nabla \times d$. There have been many interesting works on global
small solutions to the MHD system recently. See, for example, Lin-Zhang \cite{lin-zhang} and Lin-Zhang \cite{lin-zhangt}
for $n=3$ and Lin-Zhang-Xu \cite{LZX} for $n=2$.
\end{itemize}

The system \eqref{liquidcrystal} has attracted a lot of interests and generated
many interesting research works recently. Here we would like to mention a few of previous results. In dimensions two, Lin-Lin-Wang \cite{lin-lin-wang} have proved the existence of global Leray-Hopf type weak solutions to \eqref{liquidcrystal} with initial and boundary conditions, which is smooth away from finitely many possible singular times
(see Hong \cite{hong} for (\ref{liquidcrystal}) in $\Omega=\R^2$, Hong-Xin \cite{HX} and Xu-Zhang \cite{xu-zhang} for other related works). Lin-Wang \cite{lin-wang} have also proved the uniqueness for such weak solutions.  It remains a very challenging open problem to establish the existence of global Leray-Hopf type weak solutions and partial regularity of suitable weak solutions to \eqref{liquidcrystal} in dimension three. Very recently, Lin-Wang  \cite{lin-wang3} have proved the existence of global weak solutions in dimension three under the assumption that
$d_0(x)\in \mathbb S^2_+$ for a.e. $x\in\Omega$ by developing some new compactness arguments, here $
\mathbb S^2_+$ is the upper hemisphere. When $\Omega\equiv\mathbb R^3$, the local well-posedness of \eqref{liquidcrystal} was obtained for initial data $(u_0,d_0)$, where $(u_0,\nabla d_0)\in L^3_{\rm{uloc}}(\mathbb R^3)$, the space of uniformly locally $L^3$-integrable functions,  has small norms,
by Hineman-Wang \cite{hineman-wang}.  While the global well-posedness of \eqref{liquidcrystal} was obtained by Wang \cite{wang1}
for $(u_0,d_0)\in {\rm{BMO}}^{-1}\times {\rm{BMO}}$ with small norms.
A BKM type blow-up criterion was obtained
for local strong solutions to \eqref{liquidcrystal} by Huang-Wang \cite{huang-wang2} (see also Hong-Li-Xin \cite{HLX}).
More references can be found in the survey paper by Lin-Wang \cite{lin-wangs}.

It is a very interesting question whether the short time smooth solutions to the nematic liquid crystal flow
(\ref{liquidcrystal}) develop singularities in finite time. It is not hard to verify
that in dimension two, if $d$ is the heat flow of harmonic maps with finite time singularity,
constructed by \cite{chang-ding-ye}, and if we set $u\equiv 0$, then $(u, d)$ is also a solution
of (\ref{liquidcrystal}) which has a finite time singularity. However, it will
be more desirable to construct an example in which the fluid velocity field $u$ is non-trivial.
In this paper, building upon the construction by \cite{chang-ding-ye} on the heat flow of harmonic maps,
we are able to construct in dimension three the first example of solutions to (\ref{liquidcrystal})
with finite time singularity in which both the fluid field $u$ and the director field $d$ are non-trivial.
More precisely, we consider (\ref{liquidcrystal}) in the class of axisymmetric solutions without swirls and
show that for a suitably chosen domain $\Omega$ and initial-bounday data $(u_0,d_0)$, the short time
smooth solution $(u,d)$ to (\ref{liquidcrystal}) develops a finite time singularity.

Let $B_1^n\subset\mathbb R^n$ denote the unit ball centered at $0$. Now we state our first result.
\begin{theorem} \label{blowup}
Let $\Omega=B_1^2\times [0,1]$. There exists $\varphi_0\in C^\infty([0,1])$, with
$\varphi_0(0)=0$ and $|\varphi_0(1)|>\pi$, 
such that if
$$u_0({\bf x})=(x,y,-2z),$$
and
$$\displaystyle d_0({\bf x})=\Big(\frac{x}{\sqrt{x^2+y^2}}\sin\varphi_0\big(\sqrt{x^2+y^2}\big),
\frac{y}{\sqrt{x^2+y^2}}\sin\varphi_0\big(\sqrt{x^2+y^2}\big), \cos\varphi_0\big(\sqrt{x^2+y^2}\big)\Big),$$
for ${\bf x}=(x,y,z)\in\Omega$, then the short time smooth solution $(u, d, P)$ to the system (\ref{liquidcrystal})
in $\Omega$, under the initial and boundary condition:
\begin{equation}
\big(u({\bf x}, 0),d({\bf x}, 0)\big)=\big(u_0({\bf x}), d_0({\bf x})\big),\ {\bf x}\in\Omega, \label{initial}
\end{equation}
\begin{equation}\label{boundary}
\begin{cases}
\ u({\bf x},t)=u_0({\bf x}) &\ {\bf x}\in\partial\Omega, \ t>0,\\
\ d({\bf x}, t)=d_0({\bf x}) & \ {\bf x}\in \partial B_1^2\times [0,1], \ t>0,\\
\frac{\partial d}{\partial z} ({\bf x}, t)=0 & \ {\bf x}\in B_1^2\times \{0, 1\}, \ t>0.
\end{cases}
\end{equation}
must blow up at time $T_0$ for some $0<T_0=T_0(\varphi_0)<+\infty$.
\end{theorem}

Since $\Omega$ is axisymmetric and $(u_0,d_0)$ is axisymmetric without swirls, the uniqueness of short time smooth solution $(u, d)$ of (\ref{liquidcrystal}) implies that it is axisymmetric without swirls. By converting the system (\ref{liquidcrystal}) into the form of being axisymmetric without swirls, the proof of Theorem \ref{blowup} utilizes two interesting observations:
\begin{itemize}
\item[1)] The velocity $u$ is a static solution to the Navier-Stokes equation. In fact, $u(x,y,z,t)=(x,y, -2z)$
is spatial gradient of the quadratic harmonic polynomial $h(x,y,z)=\frac12x^2+\frac12y^2-z^2$.
\item [2)] The angle function $\varphi(r, t)$, associated with the orientation field $d({\bf x}, t)$, solves a drifted
version of the equation (\ref{shm}):
\begin{equation}\label{shm1}
\va_t+r\va_r=\va_{rr}+\frac{\va_r}{r}-\frac{\sin(2\va)}{2r^2}, \ (r,t)\in (0,1)\times (0,+\infty).
\end{equation}
After a suitable re-parameterization to handle the contribution from
the drifting term $r\varphi_r(r, t)$ in (\ref{shm1}), we can modify the construction
of \cite{chang-ding-ye} to build a subsolution to (\ref{shm1}) that  blows up at finite time.
This, combined with the comparison principle, yields the finite time singularity of (\ref{shm1}).
\end{itemize}
In particular, $u$ is smooth in this example. This is consistent with known results of the Navier-Stokes
equation: any local axisymmetric solution to the Navier-Stokes equation, without swirls,
is globally smooth
(see, for example, Leonardi-M$\acute{\mbox{a}}$lek-Ne$\breve{\mbox{c}}$as-Pokorn$\acute{\mbox{y}}$\ \cite{LMNP}).
Therefore, the finite time singularity arises essentially from
the orientation field $d$.

Now we would like to make a few comments related to Theorem \ref{blowup}.
\begin{remark}{\rm a) By modifying the argument of Chang-Ding \cite{chang-ding}, we will show, in Theorem 2.5 below,
that (\ref{shm1}) admits a global smooth solution if the initial-boundary value $\varphi_0\in C^\infty([0,1])$
satisfies $\varphi_0(0)=0$ and $\big\|\varphi_0\big\|_{C([0,1])}\le\pi$.\\
b) It is a natural question to ask whether the stress tensor 
$\displaystyle\mathcal S:=\nabla d\odot\nabla d-\frac12|\nabla d|^2 \mathbb I_3$ blows up in finite time, if $(u,d,P)$  is
the solution constructed by Theorem \ref{blowup}. The calculations in Section 5 seem to suggest that $\mathcal S$
may not blow up. \\
c) Since $u_0({\bf x})\cdot \nu\not=0$ on $\partial\Omega$, here $\nu$ denotes the outward unit normal of $\partial\Omega$,
it is unclear whether the solution $(u,d,P)$ constructed by Theorem \ref{blowup} enjoys the energy dissipation inequality
for $0\le t<T_0$:
\begin{equation}\label{dissipation}
\frac{d}{dt}\int_\Omega (|u|^2+|\nabla d|^2)(t)+2\int_0^t\int_\Omega \big(\mu|\nabla u|^2+\frac{\lambda}{\gamma}
|\Delta d+|\nabla d|^2 d|^2\big)\le 0.
\end{equation}
See Section 5 for more details.}
\end{remark}

It is the question from Remark 1.2 c) that motivates us to construct another example in which the solution
$(u,d, P)$ to (\ref{liquidcrystal}) develops finite time singularity, and satisfies the energy dissipation inequality (\ref{dissipation}).

In order to state it, we need some notations. Denote the north pole by ${\bf e}=(0,0,1)\in\mathbb S^2$. Set
$$C^\infty_{0,{\rm{div}}} (B^3_1,\mathbb R^3) 
:= \Big\{v \in C^\infty(B_1^3,\mathbb R^3) \ \big| \ \nabla\cdot v =  0\Big\},$$
and
$$C_{\bf e}^\infty(B_1^3,\mathbb S^2):=\Big\{d\in C^\infty(B_1^3,\mathbb S^2)\ \big| \
d={\bf e}\ {\rm{on}}\ \partial B_1^3\Big\}.$$

For continuous maps $f, g\in C\big(\overline{B_1^3},\mathbb S^2\big)$, with $f=g$ on $\partial B_1^3$,
we say that $f$ is homotopic to $g$ relative to $\partial B_1^3$ if there exists a continuous map
$\Phi\in C\big(\overline{B_1^3}\times [0,1],\mathbb S^2\big)$ such that\\
(i) $\Phi(\cdot, t)=f(\cdot)=g(\cdot)$ on $\partial B_1^3$, for all $0\le t\le 1$; and\\
(ii) $\Phi(\cdot, 0)=f(\cdot)$ and $\Phi(\cdot, 1)=g(\cdot)$ in $B_1^3$.

Now we have
\begin{theorem} \label{blow-up2} There exists 
$\epsilon_0 >0$ such that if $u_0 \in C^\infty_{0, {\rm{div}}}(B_1^3, \mathbb R^3)$ and 
$d_0 \in C_{\bf e}^\infty(B_1^3, \mathbb S^2)$ satisfies that
$d_0$ is not homotopic to the constant map ${\bf e}:B_1^3\to\mathbb S^2$ relative to $\partial B_1^3$, and
\begin{equation}\label{small_energy}
E(u_0,d_0):=\frac12\int_{B_1^3}(|u_0|^2+|\nabla d_0|^2)\le\epsilon_0^2.
\end{equation}
Then the short time smooth solution $(u, d, P) : B_1^3 \times [0, T) \to \mathbb R^3 \times\mathbb S^2$
to the nematic liquid crystal flow (\ref{liquidcrystal}),
under the initial-boundary condition
\begin{equation}\label{IBC}
\begin{cases}
(u,d)\big|_{t=0} = (u_0,d_0),\ {\rm{in}}\  B_1^3,\\
(u,d)\big|_{\partial B_1^3} = (0,{\bf e}), \ 0 < t < T,
\end{cases}
\end{equation}
 must blow up before time $T = 1$.
\end{theorem}

The following remark indicates that there are ample  examples of $(u_0,d_0)\in C^\infty_{0, {\rm{div}}}(B_1^3,\mathbb R^3)
\times C^\infty_{\bf e}(B_1^3, \mathbb S^2)$  satisfying the conditions of Theorem \ref{blow-up2}.
\begin{remark}{\rm a) Let $H(z,w)=(|z|^2-|w|^2, 2zw):\mathbb S^3\equiv \big\{(z,w)\in \mathbb C\times\mathbb C: |z|^2+|w|^2=1\big\}\to \mathbb S^2\subset \mathbb R\times \mathbb C$ be the Hopf map. Let $D_\lambda({\bf x})=\lambda {\bf x}:\mathbb R^3\to\mathbb R^3$ be the dilation map for $\lambda>0$, $\Pi:\mathbb S^3\to\overline{\mathbb R^3}$ be the 
stereographic projection map from ${\bf e}$, and $\Psi_\lambda=\Pi^{-1}\circ D_\lambda\circ \Pi:\mathbb S^3\to\mathbb S^3$. Then direct calculations imply that the Dirichlet energy of $H\circ\Psi_\lambda:\mathbb S^3\to\mathbb S^2$
satisfies
$$\lim_{\lambda\rightarrow\infty}\int_{\mathbb S^3}\big|\nabla(H\circ \Psi_\lambda)\big|^2\,d\sigma=0.$$
Moreover, it is easy to see that $H\circ\Psi_\lambda$ is not homotopic to the constant map ${\bf e}:\mathbb S^3\to\mathbb S^2$. Let $\Phi\in C^\infty\big(\overline{B_1^3}, \mathbb S^3\big)$ such that $\Phi: B_1^3\to \mathbb S^3\setminus\{\bf e\}$ is a diffeomorphism and $\Phi={\bf e}$ on $\partial B_1^3$. Now we can check that for any
$u\in C^\infty_{0,{\rm{div}}}(B_1^3,\mathbb R^3)$, since
$$\lim_{\lambda\rightarrow \infty} E(\lambda^{-1} u, H\circ\Psi_\lambda\circ\Phi)=0,$$
we can find a sufficiently large $\lambda_0>0$ depending on $u$, $H$, and $\Phi$ such that
$$(u_0,d_0):=(\lambda_0^{-1} u, H\circ\Psi_{\lambda_0}\circ\Phi):B_1^3\to\mathbb R^3\times\mathbb S^2$$
satisfies the condition (\ref{small_energy}) of Theorem \ref{blow-up2}, and $d_0$ is not homotopic to the constant map
${\bf e}$ relative to $\partial B_1^3$.\\
b) The initial data $d_0$ constructed in a) has previously been used by Ding-Wang \cite{DW} in the construction of
finite time singularity of the Landau-Lifshitz-Gilbert system modeling the continuum theory of
ferromagnetism.} 
\end{remark}
Since precise values of the parameters $\lambda$ and $\gamma$ in (\ref{liquidcrystal})
don't play a role in this paper, for simplicity we assume henceforth that
$$\lambda=\gamma=1.$$ 

The paper is organized as follows. In Section 2, we will derive the axisymmetric form of 
(\ref{liquidcrystal}), without swirls. In Section 3, we will sketch the proof of
local smooth axisymmetric solutions of (\ref{liquidcrystal1}) without swirls. 
In Section 4, we will prove the existence of global smooth axisymmetric
solutions of (\ref{liquidcrystal1}) without swirls,
when the initial data $\varphi_0$ satisfies $\varphi_0(0)=0$ and $|\varphi_0(r)|\le\pi$ for all $0\le r\le 1$.
In Section 5, we will present in details the example of solutions to (\ref{liquidcrystal1})
with finite time singularity for suitably chosen initial data $\phi_0\in C^\infty([0,1])$, 
with $\varphi_0(0)=0$ and $|\varphi_0(1)|>\pi$. In Section 6, we will outline the proof of Theorem
\ref{blow-up2}.

\section {Axisymmetric form of (\ref{liquidcrystal}) without swirls}
\setcounter{equation}{0}
\setcounter{theorem}{0}
In this section, we will derive the axisymmetric form of (\ref{liquidcrystal}) without
swirls. We would like to mention that Dong-Lei \cite{dong-lei} have constructed 
a global smooth axisymmetric solution of (\ref{liquidcrystal}) in dimension two.

Let $(r,\theta, z)$ denote the cylindrical
coordinates of $\mathbb R^3$, and set
$$\e^r=(\cos\theta,\sin\theta,0),\quad \e^{\theta}=(-\sin\theta,\cos\theta,0),\quad \e^3=(0,0,1)$$
as the canonical orthonormal base of $\mathbb R^3$ in the cylindrical coordinates.
For $\alpha\in [0,2\pi]$, let $R_\alpha\in {\rm{SO}}(3)$ denote the rotation map 
of angle $\alpha$ with respect to the $z$-axis. 
Recall that a vector field $v:\mathbb R^3\to\mathbb R^3$ is axisymmetric if
$$R_\alpha^{-1}\circ v\circ R_\alpha =v, \ \forall\ \alpha\in [0,2\pi].$$
Hence any axisymmetric vector field $v$ can be written as
\beq\label{axis1.1}
v(r,\theta,z)=v^r(r, z)\e^r+v^{\theta}(r,z)\e^{\theta}+v^3(r,z)\e^3.
\eeq
If, in addition, $v^\theta\equiv 0$, we say $v$ is axisymmetric without swirls.

A solution $(u,d, P)$ of the nematic liquid crystal flow equation (\ref{liquidcrystal}) is said to
be axisymmetric without swirls, if 
$$\begin{cases} u(r,\theta, z,t)=u^{r}(r,z,t)\e^r+u^{3}(r,z,t)\e^3,\\
d(r,\theta, z,t)=\sin\va(r,z,t) \e^r+\cos\va(r,z,t) \e^3,\\
P(r,\theta,z,t)=P(r,z,t).
\end{cases}
$$
A domain $\Omega\subset\mathbb R^3$ is axisymmetric if it is invariant under a rotation
map $R_\alpha$ for any $\alpha\in [0,2\pi]$. Now we have
\begin{lemma} For any axisymmetric domain $\Omega\subset\mathbb R^3$,
if $(u,d,P)$ is an axisymmetric without swirl solution of the system (\ref{liquidcrystal}) in $\Omega\times \mathbb R_+$,
then $(u^r, u^3, \varphi, P)$ solves 
\begin{equation}\label{liquidcrystal1}\displaystyle
\begin{cases}
\displaystyle\frac{\D u^r}{Dt}-\mu\del u^r+\frac{\mu}{r^2}u^r+P_r
=-\big(\del\va-\frac{\sin(2\va)}{2r^2}\big)\va_r,\\
\ \ \ \ \ \ \ \ \ \displaystyle\frac{\D u^3}{Dt}-\mu\del u^3+P_z
=-\big(\del\va-\frac{\sin(2\va)}{2r^2}\big)\va_z,\\
\qquad\qquad\ \displaystyle\frac{1}{r}(ru^r)_r+(u^3)_z=0,\\
\qquad\qquad\qquad\ \ \displaystyle\frac{\D \va}{Dt}-\del\va
=-\frac{\sin(2\va)}{2r^2},
\end{cases}
\end{equation}
where
\bex
\begin{cases}
\displaystyle\frac{\D}{Dt}:=\pa_t+u^r\pa_r+u^3\pa_z,\\
\ \ \displaystyle\del:=\pa^2_r+\frac{1}{r}\pa_r+\pa^2_z.
\end{cases}
\eex
\end{lemma}

\pf Let's first show (\ref{liquidcrystal1})$_3$. By (\ref{liquidcrystal})$_2$ and the definition of $u$, we have
\begin{eqnarray*}
0&=&\pa_{x}\left(\frac{u^rx}{r}\right)+\pa_{y}\left(\frac{u^ry}{r}\right)+\pa_{z}\left(u^3\right)\\
&=&\frac{\pa_{x}(u^rx)+\pa_{y}(u^ry)}{r}-\frac{u^r(x\pa_{x}r+y\pa_{y}r)}{r^2}+\pa_{z}\left(u^3\right)\\
&=&\frac{2u^r}{r}+\pa_ru^r+u^3_z-\frac{u^r}{r}=\frac{(ru^r)_r}{r}+u^3_z.
\end{eqnarray*}
For \eqref{liquidcrystal1}$_4$, since $d$ can be written as
\bex
d=(\cos\theta\sin\va,\sin\theta\sin\va,\cos\va),
\eex
direct calculations imply that
\bex
\begin{cases}
d_t=(\cos\theta\cos\va,\sin\theta\cos\va,-\sin\va) \va_t,\\
d_r=(\cos\theta\cos\va,\sin\theta\cos\va,-\sin\va)\va_r,\\
d_{\theta}=(-\sin\theta\sin\va,\cos\theta\sin\va,0),\\
d_z=(\cos\theta\cos\va,\sin\theta\cos\va,-\sin\va)\va_z,
\end{cases}
\eex
and
$$
\begin{cases}
d_{rr}=\left((\cos\va\,\va_{rr}-\sin\va \va_r^2)(\cos\theta,\sin\theta),
-(\sin\va\,\va_{rr}+\cos\va\va_r^2)\right),\\
d_{\theta\theta}=(-\cos\theta\sin\va,-\sin\theta\sin\va,0),\\
d_{zz}=\left((\cos\va\,\va_{zz}-\sin\va \va_z^2)(\cos\theta,\sin\theta),
-(\sin\va\,\va_{zz}+\cos\va \va_{z}^2)\right).
\end{cases}
$$
Hence we have
\bex
|\nabla d|^2=\va_r^2+\frac{1}{r^2}\sin^2\va+\va_z^2,
\eex
and
\bex
\Delta d+|\nabla d|^2d=\big(\va_{rr}+\va_{zz}+\frac{1}{r}\va_r-\frac{\sin(2\va)}{2r^2}\big)
\big(\cos\theta\cos\va,\sin\theta\cos\va,-\sin\va\big).
\eex
By the definition of $u$,  we also have
\bex
u\cdot\nabla d=(u^r\pa_r+u^3\pa_z)d
=\left(u^r\va_r+u^3\va_z\right)
(\cos\theta\cos\va,\sin\theta\cos\va,-\sin\va).
\eex
Putting these identities into (\ref{liquidcrystal})$_3$, we obtain
\bex
\va_t+u^r\va_r+u^3\va_z=\big(\va_{rr}+\va_{zz}+\frac{1}{r}\va_r-\frac{\sin(2\va)}{2r^2}\big).
\eex
This yields the equation (\ref{liquidcrystal1})$_4$.
Now we want to derive the momentum equations (\ref{liquidcrystal1})$_{1,2}$.
Since
\beq\label{eqn_er}
\Delta \e^r=-\frac{1}{r^2}\e^r,
\eeq
it follows from (\ref{liquidcrystal})$_1$ and (\ref{eqn_er}) that
\beq\label{axis1.13}
\begin{split}
\big(\frac{\D u^r}{Dt}-\mu\del u^r+\frac{\mu}{r^2}u^r+P_r\big)\e^r
+\big(\frac{\D u^3}{Dt}-\mu\del u^3+P_z\big)\e^3=-\Delta d\cdot\nabla d.
\end{split}
\eeq
From (\ref{eqn_er}), we also have 
$$
\Delta d=\big(\Delta-\frac{1}{r^2}\big)(\sin\va) \e^r+\Delta (\cos\va)\e^3.
$$
Since
\bex
\nabla d=\nabla(\sin\va) \e^r+\nabla(\cos\va)\e^3+\sin\va\nabla\e^r,
\eex
and
\bex
\langle\nabla \e^r,\e^r\rangle=\langle\nabla \e^r, \e^3\rangle=0,
\eex
we have
\beq\label{axis1.16}
\begin{split}
-\Delta d\cdot\nabla d=&-\big(\Delta-\frac{1}{r^2}\big)(\sin\va)\nabla(\sin\va)
-\Delta (\cos\va)\nabla(\cos\va)\\
=&\big(\Delta (\cos\va)\sin\va-\big(\Delta-\frac{1}{r^2}\big)(\sin\va)\cos\va\big)
(\va_r \e^r+\va_z \e^3)\\
=&-\big(\va_{rr}+\va_{zz}+\frac{1}{r}\va_r-\frac{\sin(2\va)}{2r^2}\big)(\va_r \e^r+\va_z \e^3).
\end{split}
\eeq
Putting (\ref{axis1.16}) into (\ref{axis1.13}) yields (\ref{liquidcrystal})$_{1,2}$.
\endpf

\smallskip
In order to construct a solution of (\ref{liquidcrystal1}) with finite time singularity,
we further consider the domain to be the round cylinder $\Omega=B_1^2\times [0,1]$ 
and an axisymmetric, without swirl, solution $(u,P, d)$ in the special form:
$$\begin{cases}
u(r,\theta, z,t):=v(r,t)\e^r+w(z,t)\e^3,\\
d(r,\theta, z, t):=\sin\va(r,t) \e^r+\cos\va(r,t) \e^3,\\
P(r,\theta, z,t):=Q(r,t)+R(z,t).
\end{cases}
$$
Then \eqref{liquidcrystal1} becomes
\begin{equation}\label{liquidcrystal2}
\begin{cases}
\displaystyle v_t+vv_r-\mu\big( v_{rr}+\frac{v_r}{r}-\frac{1}{r^2}v\big)+Q_r
=\big(\va_{rr}+\frac{\va_r}{r}-\frac{\sin(2\va)}{2r^2}\big)\va_r, \  r\in [0,1], \\
\qquad\qquad\quad \ w_t+ww_z-\nu w_{zz}+R_z
=0, \ \qquad\qquad \ \ \ \ \ \ \ \ \ \ \ \ \ \ \ \ \ \ z\in [0,1], \\
\qquad\qquad\qquad\qquad\ \ \ \ \ \ \displaystyle\frac{1}{r}(rv)_r+w_z=0,
\ \ \ \ \ \ \ \ \ \ \ \ \ \ \ \    (r,z)\in [0,1]\times [0,1],\\
\qquad\qquad\ \ \ \ \ \ \ \ \ \ \ \ \ \ \ \ \ \ \ \ \ \  \displaystyle\va_t+v\va_r=\va_{rr}+\frac{\va_r}{r}-\frac{\sin(2\va)}{2r^2},
 \qquad\ r\in [0,1].
\end{cases}
\end{equation}
The initial condition for $(v,w,\phi)$ reduces to
\begin{equation}\label{sbinitial}
\begin{cases}
v|_{t=0}=r,\ \ \ \ \ \ \ 0\le r\le 1,\\
w|_{t=0}=-2z, \ \ 0\le z\le 1,\\
\va|_{t=0}=\va_0(r),\ 0\le r\le 1,
\end{cases}
\end{equation}
for some $\va_0\in C^\infty([0,1])$, with $\va_0(0)=0$.
The boundary condition is
\beq\label{sboundary}
\begin{split}
\begin{cases} v(0,t)=0\\
v(1,t)=1,
\end{cases} & \quad \begin{cases} w(0,t)=0\\ w(1,t)=-2,\end{cases}
\quad \begin{cases} \va(0,t)=0\\ \va(1,t)=\va_0(1).\end{cases}
\end{split}
\eeq
\begin{remark}{\rm
It is easy to verify that  the initial and boundary condition (\ref{sbinitial})-\eqref{sboundary} for
$(v, w, \phi)$ is equivalent to  the initial condition (\ref{initial})
and the boundary condition (\ref{boundary}) for $(u,d)$.}
\end{remark}

Now, we are ready to state our main results on the system
(\ref{liquidcrystal2}), under the initial and boundary condition
(\ref{sbinitial}) and (\ref{sboundary}). 

The first result asserts that for any short time smooth solution $(v,w,\varphi, Q,R)$ to
(\ref{liquidcrystal2}) along with the initial-boundary condition (\ref{sbinitial}) and (\ref{sboundary}), $(v,w)$ is static
and $Q, R$ are also static (up to a time-dependent constant). In fact, we have
\begin{lemma} \label{static} For $0<T\le +\infty$, suppose that $v, w\in C^\infty([0,1]\times [0, T))$ satisfies
\begin{equation}\label{div}
\displaystyle\frac{1}{r}(rv)_r+w_z=0, \ (r,z)\in [0,1]\times [0,1],
\end{equation}
and
\beq\label{2pt}
\begin{split}
\begin{cases} v(0,t)=0\\
v(1,t)=1,
\end{cases} & \quad \begin{cases} w(0,t)=0\\ w(1,t)=-2.\end{cases}
\end{split}
\eeq
Then $v(r,t)=r$ for any $(r,t)\in [0,1]\times [0,T)$, and $w(z,t)=-2z$ for any $(z,t)\in [0,1]\times [0,T)$.
\end{lemma} 
\pf  Differentiating \eqref{div} with respect to $z$ yields
\bex
w_{zz}(z,t)=0
\eex
so that $w(z,t)=a_1(t)z+a_2(t)$ for some functions $a_1(t)$ and $a_2(t)$. 
Since $w(0,t)=0$ and $w(1,t)=-2$, we see that
$a_2(t)\equiv 0$ and $a_1(t)\equiv -2$. Thus $w(z,t)=-2z.$

Similarly, differentiating \eqref{div} with respect to $r$ yields
\bex
\big(\frac{1}{r}(rv)_r\big)_r(r,t)=0,
\eex
this implies that $rv(r,t)=b_1(t)r^2+b_2(t)$ for some functions $b_1(t)$ and $b_2(t)$. 
Since $v(0,t)=0$ and $v(1,t)=1$, we  see that
$b_2(t)\equiv 0$ and $b_1(t)\equiv 1$. Thus $v(r,t)=r.$
The proof is complete.
\endpf\\

Next we have 
\begin{lemma} \label{drift-hm0} For $0<T\le +\infty$, $(v,w, \varphi, Q, R)\in C^\infty([0,1]\times [0,T))$ solves
(\ref{liquidcrystal2}), under (\ref{sbinitial}) and (\ref{sboundary}) iff
$v(r,t)=r$, $w(z,t)=-2z$,  
\begin{equation}\label{pressure}
\begin{cases}
R(z,t)=-2z^2+c_1(t)\\
\displaystyle Q(r,t)
=-\int_0^r\big(\va_{rr}+\frac{\va_r}{r}-\frac{\sin(2\va)}{2r^2}\big)\va_r\,dr-\frac{r^2}2+c_2(t),
\end{cases}
\end{equation}
for some $c_1, c_2\in C^\infty([0,T))$, and
\begin{equation}\label{drift-hm}
\begin{cases}
\displaystyle\va_t+r\va_r=\va_{rr}+\frac{\va_r}{r}-\frac{\sin(2\va)}{2r^2}, \ 0<r<1, \\
\varphi(r,0)=\varphi_0(r), \ 0<r<1,\\
\varphi(0,t)=0, \ \ \varphi(1,t)=\varphi_0(1), \ t>0.
\end{cases}
\end{equation}
\end{lemma}
\pf Applying lemma \ref{static} to (\ref{liquidcrystal2})$_2$ yields that 
$$R_z+4z=0,$$
this, by integration, implies 
\bex
R(z,t)=-2z^2+c_1(t),
\eex 
for some function $c_1\in C^\infty([0,T))$. 
Apply lemma \ref{static} to (\ref{liquidcrystal2})$_{1,4}$, we obtain
\begin{align}
Q_r
=&-\big(\va_{rr}+\frac{\va_r}{r}-\frac{\sin(2\va)}{2r^2}\big)\va_r-r,  \label{eqQ}\\
\va_t+r\va_r=&~\va_{rr}+\frac{\va_r}{r}-\frac{\sin(2\va)}{2r^2}.   \label{eqva}
\end{align}
Integrating (\ref{eqQ}) yields (\ref{pressure})$_1$, while (\ref{pressure})$_2$ follows from 
(\ref{eqva}), (\ref{sbinitial}),  and (\ref{sboundary}). 
\endpf\\

It is readily seen that whether (\ref{liquidcrystal2}), under (\ref{sbinitial}) and (\ref{sboundary}),
admits a global smooth solution is equivalent to whether (\ref{drift-hm}) admits a global smooth
solution. For this, we are able to extend Chang-Ding \cite{chang-ding} and Chang-Ding-Ye \cite{chang-ding-ye} on
the heat flow of harmonic maps and obtain the following two results.

\btm\label{th1}{\it Suppose $\phi_0\in C^\infty([0,1])$
satisfies $\va_0(0)=0$ and $|\va_0(r)|\le \pi$ for all $r\in [0,1]$.
Then there is a unique, global smooth solution $(v,w,\va)$
to \eqref{liquidcrystal2}, under \eqref{sbinitial} and \eqref{sboundary}.}
\etm

\btm\label{th3}{\it There exists $\phi_0\in C^\infty([0,1])$, with $\phi_0(0)=0$
and $|\va_0(1)|> \pi$, such that
the short time smooth solution $(v,w,\va)$ to \eqref{liquidcrystal2},
under \eqref{sbinitial} and \eqref{sboundary},
must blow up at $T_0$ for some
$0< T_0=T_0(\phi_0)<+\infty$. More precisely,  $\va_r(0,t)\rightarrow\infty$ as $t\rightarrow T_0^-$.
}
\etm

\begin{remark}
{\rm It is readily seen that Theorem \ref{blowup} follows directly from Theorem \ref{th3}.
Since $\va(r,  t)$ blows up at $r=0$ and $t=T_0$, we see that
$d(r, \theta, z, t)$ blows up at $r=0$, $z\in [0,1]$, and $t=T_0$. Thus  the singular set of $d$
at the first singular time $T_0$ includes the line segment $\{0\}\times [0,1]\subset \Omega$.}
\end{remark}

The proofs of Theorem 2.5 and Theorem 2.6 will be given in Section 4 and Section 5 respectively.

\section {Existence of short time smooth solutions}
\setcounter{equation}{0}
\setcounter{theorem}{0}

In this section we will establish the existence of short time smooth solutions
to the system \eqref{liquidcrystal2}, under the initial and boundary conditions
(\ref{sbinitial}) and (\ref{sboundary}). By lemma \ref{drift-hm0}, it suffices to
prove the local solvability of the drifted-harmonic map equation (\ref{drift-hm}).

\btm[{\rm local existence}]\label{locexth}
{\it For any $\varphi_0\in C^\infty([0,1])$, with $\varphi_0(0)=0$,
there exist $t_0>0$ and a unique smooth solution
 $\va(r,t)$ on $[0,1]\times [0,t_0)$ to the equation (\ref{drift-hm}).}
\etm

\pf First observe that if $\varphi$ solves  (\ref{drift-hm}),
then
$$d(x,y,t):=\Big(\sin\va\big(\sqrt{x^2+y^2},t\big)\frac{x}{\sqrt{x^2+y^2}},
\sin\va\big(\sqrt{x^2+y^2},t\big)\frac{y}{\sqrt{x^2+y^2}},
\cos\va\big(\sqrt{x^2+y^2},t\big)\Big),$$
for $(x,y)\in B_1^2$ and $t\ge 0$, solves
\beq\label{eqd2}
d_t+(x,y)\cdot\nabla d=\Delta d+|\nabla d|^2 d,\ (x,y)\in B_1^2,\ t>0,
\eeq
with the initial condition
\beq\label{d2initial}
d\big|_{t=0}=\Big(\sin\va_0\big(\sqrt{x^2+y^2}\big)\frac{x}{\sqrt{x^2+y^2}},
\sin\va_0\big(\sqrt{x^2+y^2}\big)\frac{y}{\sqrt{x^2+y^2}}, \cos\va_0\big(\sqrt{x^2+y^2}\big)\Big),
\eeq
and the boundary condition
\beq\label{d2boundary}
d\big|_{\partial B_1^2}=\big(\sin \va_0(1) x, \sin\va_0(1)y, \cos \va_0(1) \big).
\eeq
By the standard contraction mapping theorem, we can prove that there exist $t_0>0$
and a unique smooth solution $d\in C^\infty(B_1^2\times [0,t_0), \mathbb S^2)$
to \eqref{eqd2}, along with \eqref{d2initial} and \eqref{d2boundary}
(see, e.g., \cite{lin-lin-wang} Section 3 for a detailed proof).

Now we need to argue that such a solution $d(x,y,t)$ to (\ref{eqd2})-(\ref{d2initial})-(\ref{d2boundary})
is axisymmetric, i.e.,
\begin{equation}\label{axisym}
d(x, y, t)(:=d(r,\theta, t))=\big(\sin \va(r,t)\cos\theta, \sin\va(r,t)\sin\theta, \cos \va(r,t)\big).
\end{equation}
Assume (\ref{axisym}) is true. Then a straightforward calculation as in Section 2 implies that
$\varphi$ solves (\ref{drift-hm}).
To see (\ref{axisym}), first note that we can always write
\beq\label{loc1}
d(x,y,t)=
\big(\sin \va(r,\theta,t)\cos\chi(r,\theta,t),\, \sin\va(r,\theta,t)\sin\chi(r,\theta,t),\, \cos \va(r,\theta,t) \big)
\eeq
for a pair of functions $\va$ and $\chi$.  For $\alpha\in [0, 2\pi]$, let $R_\alpha\in {\rm{SO}}(3)$ 
be the rotation map in $\R^3$ of angle $\alpha$ with respect to the $z$-axis, and $S_\alpha\in {\rm{SO}}(2)$
be the rotation  map in $\R^2$ of angle $\alpha$ with respect to $0$.
Since the initial and boundary values $d_0$ of $d$ are axisymmetric, i.e., for any $\alpha\in [0,2\pi]$,
$$d_0(S_\alpha(x,y))=R_\alpha \circ d_0(x,y), \ (x,y)\in \overline{B_1^2},$$
we can check that $
R^{-1}_{\alpha}\circ d\big(S_\alpha(x,y),t\big), \ \forall \alpha\in [0,2\pi],$
also solves the initial and boundary value problem of the drifted harmonic map equation (\ref{drift-hm}). 
By the uniqueness of solutions to (\ref{drift-hm}), this implies
that 
$$
d(x,y, t)=R^{-1}_{\alpha}\circ d\big(S_\alpha(x,y),t\big), \ \forall\alpha\in [0,2\pi], \ (x,y)\in B_1^2,
$$
or equivalently,
\beq\label{loc2}
d(x,y,t)=
\big(\sin \va(r,t)\cos(\eta(r,t)+\theta),\, \sin\va(r,t)\sin(\eta(r,t)+\theta),\, \cos \va(r,t) \big)
\eeq
for some functions $\varphi(r,t)$ and $\eta(r,t)$. Direct calculations, similar to Section 2,
imply that $\va$ and $\eta$ satisfy 
\beq\label{loc3}
\begin{cases}
\displaystyle
\va_t+r\va_r=\big(\va_{rr}+\frac{\va_r}{r}\big)-\big(\eta_r^2+\frac{1}{r^2}\big)\sin\va\cos\va,\\
\displaystyle\sin\va(\eta_t+r\eta_r)=\sin\va\big(\eta_{rr}+\frac{\eta_r}{r}\big)+2\va_r\eta_r\cos\va.
\end{cases}
\eeq
Now we need to show

\bigskip
\noindent{\bf Claim 1}. {\it $\eta(r,t)=0$ for $(r,t)\in[0,1]\times[0,t_0)$.}
\bigskip
\\
Note that when $\sin\va(r,t)=0$, we can simply define $\eta(r,t)=0$ because
any value of $\eta(r,t)$ gives the same value of $d$ by \eqref{loc2}. Set
$$
U=\Big\{(r,t)\in[0,1]\times[0,t_0)\ \big|\  \sin\va(r,t)\neq0\Big\},
$$
and define
$$
h(r,t)=\eta(r,t)\sin\va(r,t), \ \ (r,t)\in U.
$$
Then by \eqref{loc3}, $h(r,t)$ satisfies
\beq\label{loc4}
h_t+rh_r=h_{rr}+\frac{h_r}{r}+\big[\va_r^2-\big(\eta_r^2+\frac{1}{r^2}\big)\cos^2\va\big]h,
\eeq
under the initial and boundary conditions:
\beq\label{loc5}
h(r,t)=0, \ {\rm{on}}\ \partial U\cap \{0\le t<t_0\}.
\eeq
Since $\va_r$ is bounded in $U$, replacing $h$ by $h\,e^{-ct}$ for a sufficiently large $c>0$ we may
assume the coefficient of $h$ in \eqref{loc4} is negative. Therefore by the maximum principle,
we conclude that
$h\equiv 0$ and consequently $\eta\equiv0$ on $U$ (see \cite{chang-ding} Lemma 2.2, or \cite{friedman} and \cite{pw}).  Hence Claim 1 holds and the proof of Theorem 3.1 is complete.
\endpf


\section {Existence of global smooth solutions and proof of Theorem \ref{th1}}
\setcounter{equation}{0}
\setcounter{theorem}{0}

This section is devoted to the proof of Theorem \ref{th1} on the existence of global smooth
solutions when the initial data $\phi_0$ satisfies $\phi_0(0)=0$ and $|\phi_0(r)|\le\pi$ for all $r\in [0,1]$. 
The proof is motivated by \cite{chang-ding}.

We start with

\blm[{\rm maximum principle}]\label{maxp}
{\it For any smooth solution $\va(r,t):[0,1]\times [0,T)\rightarrow \mathbb R$ to the drifted harmonic
map equation \eqref{drift-hm}, if
$$
\va_0(0)=0,  \  -\pi\leq \va_0(r)\leq \pi, \ \forall \ r\in [0,1],
$$
then
\beq
-\pi< \va(r,t)<\pi, \ \forall\ r\in (0,1)\ {\rm{and}}\ t\in (0,T).
\eeq
}
\elm

\pf We only present the proof for $\va<\pi$, since the other half of the inequality (4.1)
can be proved similarly. Set $\tilde \va=\pi-\va$. Then $\tilde\va$ satisfies
\beq\label{lmpf1}
\displaystyle\tilde\va_t+r\tilde\va_r=~\tilde\va_{rr}+\frac{\tilde\va_r}{r}+p_1(r,t)\tilde\va,
\eeq
where
$$
p_1(r,t):=\frac{\sin(2\tilde\va)}{2\tilde\va r^2},
$$
and 
$$
\tilde\va(r,0)\geq 0, ~\tilde\va(1,t)\geq 0 \mbox{ and } \tilde\va(0,t)=\pi.
$$
Thus for any $t\in(0,T)$, there exists a small $r_1\in(0,1)$ such that $\tilde\va(r,\tau)>0$ for any $(r,\tau)\in(0,r_1]\times(0,t]$. On $(r_1,1)\times(0,t]$, by the fact that $p_1(r,\tau)$ is bounded, we can prove $\tilde\va>0$ by the standard maximum principle (see \cite{friedman} or \cite{pw}). Therefore, we conclude that $\tilde\va>0$
or equivalently $\va<\pi$ on $(0,1)\times(0,T)$.
\endpf

\blm[{\rm comparison principle}]\label{comp}
{\it Suppose the functions $\va$, $f$ and $g$ are smooth solution, subsolution and supersolution to \eqref{eqva} on $[0,1]\times[0,T)$ respectively, and
$$f(r,t)\le \va(r,t)\le g(r,t) \ {\rm{on}}\ ([0,1]\times\{0\})\cup (\{0,1\}\times (0,T)).$$
Then we have
\beq
f(r,t)\le \va(r,t) \le g(r,t), \ \forall \ (r,t)\in [0,1]\times[0,T).
\eeq
}
\elm

\pf  Set $\bar f=f-\va$. Then $\bar f$ satisfies
\beq\label{lmpf2}
\displaystyle
\bar f_t+r\bar f_r\ge ~\bar f_{rr}+\frac{\bar f_r}{r}+p_2(r,t)\bar f,
\eeq
where
$$
\displaystyle p_2(r,t):=-\frac{\sin(2f(r,t))-\sin(2\va(r,t))}{2 r^2(f(r,t)-\va(r,t))},
$$
and
$$
\bar f(r,t)\le 0 \ {\rm{on}}\ ([0,1]\times\{0\})\cup (\{0,1\}\times (0,T)).
$$
For any $t\in(0,T)$, there exists a small $r_2\in(0,1)$ such that $p_2(r,\tau)<0$ on $(r,\tau)\in(0,r_2)\times(0,t)$. Combining with the fact $p_2(r,\tau)$ is bounded on $(r_2, 1)\times(0,t)$, we conclude that $p_2(r,\tau)$ is bounded from above on $(r,\tau)\in(0,1)\times(0,t)$. By the standard maximum principle (see \cite{friedman} or \cite{pw}), 
we conclude that $\bar f\le 0$ or $\va\ge f$ on $[0,1]\times(0,T)$. 
Similarly, one can prove $\phi\le g$.
\endpf

\bigskip
\noindent{\bf Proof of Theorem \ref{th1}.}\quad 
To prove the existence of global smooth solutions of (\ref{drift-hm}),
we need to construct suitable supersolutions and subsolutions to \eqref{eqva}. Denote
\beq\label{pf1}
\overline\va(r,c)=2\arctan \left(\frac{r}{c}\right)
\quad\mbox{and}\quad
\underline\va(r,c)=2\arctan \left(-\frac{r}{c}\right),
\eeq
for some positive constant $c$. It is easy to see that $\overline\va(r,c)$ and $\underline\va(r,c)$ are smooth functions 
in $[0,1]$. Direct calculations  give
\beq\label{pf2}
\overline\va_t+r\overline\va_r-\overline\va_{rr}-\frac{\overline\va_r}{r}-\frac{\sin(2\overline\va)}{2r^2}
=\frac{2rc}{c^2+r^2}\geq 0, \ r\in [0,1].
\eeq
Thus $\overline{\va}$ is a supersolution of (\ref{eqva}).
Denote
$$\eta_0=\pi-\max\limits_{0\leq r\leq 1}|\va_0(r)|.$$
By the assumption on $\varphi_0$, we have $\eta_0\ge 0$ and $|\va_0|\leq\pi-\eta_0$.
Since $\overline\va_r(0,c)=\frac{2}{c}$ and $\overline\va(0,c)=\va_0(0)=0$, we can find a sufficiently small $c>0$ 
such that
\beq\label{pf3}
\overline\va(r,c)\geq\va_0(r)
\eeq
 for any $r\in[0,1]$, with equality iff $r=0$. Similarly, we can prove that $\underline\va(r,c)$ is a subsolution to \eqref{eqva}
 and 
 \beq\label{pf3.1} \varphi_0(r)\ge \underline\va(r,c) \ {\rm{in}} [0,1]
 \eeq
for a sufficiently small $c>0$. By lemma \ref{comp}, 
 we can conclude that $\overline\va(r,c)\ge\va(r,t)\ge\underline\va(r,c)$ for $r\in [0,1]$ and $t>0$.

Suppose $T>0$ is the maximum time interval for $\va$. For any $r_0\in (0,1)$ and $t_1\in(0,T)$, by
the standard regularity theory of parabolic equations, we can prove
\beq\label{pf4}
\big\|\va\big\|_{C^k([r_0,1]\times(t_1,T))}\leq Cr_0^{-k}, \ \forall k\ge 1.
\eeq
This implies that the possible singularity of the  solution $\phi$ can only happen at $r=0$.
Suppose $T<+\infty$. Then
$$d(x,y,t)=(\sin\va\cos\theta,\sin\va\sin\theta,\cos\va)$$
blows up at $\big((0,0),T\big)$. By the standard blowing up argument (cf. \cite{lin-lin-wang} Theorem 1.3),
there exist $(x_m,y_m)\rightarrow (0,0)$, $t_m\uparrow T$, $r_m\downarrow 0$ as $m\rightarrow+\infty$ and a nonconstant smooth harmonic map $\omega:\mathbb R^2\rightarrow \mathbb S^2$ such that
\beq\label{pf5}
d_m(x,y,t):=d(r_mx,r_my, t_m+r_m^2t)\rightarrow \omega\quad \mbox{in } C^2_{\mbox{loc}}(\mathbb R^2\times(-\infty,0]).
\eeq
Since, by lemma \ref{maxp},
$$-\pi< \va(r,t)<\pi, \ r\in [0, 1], \ t_m\le t<T,
$$
we have that the third component, $\omega^3$, of $\omega$ satisfies
$$
-1<\omega^3(x,y)< 1\quad \mbox{for any }(x,y)\in \mathbb R^2.
$$
This implies that  $\omega\in C^\infty(\mathbb R^2,\mathbb S^2)$ is a nontrivial harmonic map,
with finite energy and degree zero, which is impossible.
Therefore, we conclude that $T=+\infty$. This proves Theorem \ref{th1}.
\qed

\section {Finite time singularity and proof of Theorem \ref{th3}}
\setcounter{equation}{0}
\setcounter{theorem}{0}

In this section, we will modify the construction by \cite{chang-ding-ye} to show
the existence of solutions with finite time singularity of (\ref{drift-hm}) for some suitably chosen
initial data $\varphi_0$ with $\varphi_0(0)=0$ and $|\varphi_0(1)|>\pi$.\\

\noindent{\bf Proof of Theorem \ref{th3}}\quad
We adopt some ideas of \cite{chang-ding-ye} to construct suitable barrier functions from below.
Without loss of generality, we assume that $\va_0(1)>\pi$.
Denote the scalar-valued tension field by
$$\displaystyle
\tau(\phi):=\phi_{rr}+\frac{1}{r}\phi_r-\frac{\sin \phi\cos \phi}{r^2}.
$$
There are a family of smooth solutions of  $\tau(\phi)=0$,
 with $\phi(0)=0$ and $\phi(r)>0$, given by
\beq\label{axis6.1}
\phi(r,\beta)=\arccos \left(\frac{\beta^2-r^2}{\beta^2+r^2}\right)
\quad \mbox{or} \quad
\phi(r,\beta)=2\arctan \left(\frac{r}{\beta}\right)
\eeq
for $\beta>0$.

In order to  handle the drift term $r\va_r$ in the equation \eqref{eqva},
we need to reparameterize $\beta$ by $e^{t}\beta$ in \eqref{axis6.1}:
\beq\label{axis6.2}
\phi(r,\beta,t)=\arccos \left(\frac{e^{2t}\beta^2-r^2}{e^{2t}\beta^2+r^2}\right)
\quad \mbox{or} \quad
\phi(r,\beta,t)=2\arctan \left(\frac{r}{e^t\beta}\right),
\eeq
which is also a solution of $\tau(\phi)=0$ for any $t\in[0,+\infty)$.
For any $\ve\in (0,1)$ and $\mu>0$, let $a=1+\ve$ and
\beq\label{axis6.3}
\theta(r,\mu,t)=2\arctan \left(\frac{r^a}{e^{at}\mu}\right).
\eeq
Then $\theta(r,\mu,t)$ satisfies
\beq\label{axis6.4}\displaystyle
\theta_{rr}+\frac{1}{r}\theta_r-\frac{a^2\sin\theta\cos\theta}{r^2}=0.
\eeq
Choosing $\mu$ large enough so that $\theta(r,\mu,t)$ is small enough and
\beq\label{axis6.5}
\cos \theta(r,\mu,t)\geq \frac{1}{1+\ve}
\eeq
for any $r\in[0,1]$ and $t\in[0,+\infty)$.
We will  look for a subsolution $f$ of \eqref{eqva} in the form:
\beq\label{axis6.6}
f(r,t)=\phi(r,\beta(t),t)+ \theta(r,\mu,t),
\eeq
where $\beta(t)$ solves
\beq\label{axis6.7}
\begin{cases}
\frac{d\beta}{dt}=-\delta e^{-2t}\beta^{\ve},\\
\beta(0)=\beta_0,
\end{cases}
\eeq
where $\delta$ and $\beta_0$ are positive constants to be determined later.

\medskip
\noindent{\bf Claim 2}. {\em If $0<2\beta_0^{1-\ve}<\delta(1-\ve)$, then there exists $0<T_0<+\infty$ such that
\beq\label{sub_blowup}
\lim_{t\uparrow T_0^{-}}f_r(0,t)=+\infty.
\eeq}
To see (\ref{sub_blowup}), we solve the ordinary differential equation \eqref{axis6.7}
and  obtain
\beq\label{axis6.14}
\beta^{1-\ve}=\frac{\delta(1-\ve)}{2}(e^{-2t}-1)+\beta_0^{1-\ve}.
\eeq
Set
$$
T_0:=\frac{1}{2}\ln \left(\frac{\delta(1-\ve)}{\delta(1-\ve)-2\beta_0^{1-\ve}}\right)>0.
$$
Then 
$$\beta(t)\rightarrow 0\ {\rm{as}}\ t\rightarrow T_0^{-}.$$
By \eqref{axis6.10} below, we have
$$
f_r(r,t)=\frac{2\beta e^t}{e^{2t}\beta^2+r^2}+\frac{2a\mu r^{\ve}e^{at}}{\mu^2e^{2at}+r^{2a}}
$$
so that
 $$
f_r(0,t)=\frac{2}{e^{t}\beta(t)}\rightarrow +\infty,\quad \mbox{as } t\rightarrow T_0^-.
$$

\bigskip
\noindent{\bf Claim 3}. {\it There exists $\delta>0$ such that
$f$  is a subsolution of (\ref{drift-hm}), i.e., $f$ satisfies
 $$f_t+rf_r\leq \tau(f).$$}

\noindent In fact, by \eqref{axis6.1} and \eqref{axis6.4}, we have
\beq\label{axis6.8}
\begin{split}
\tau(f)=&\frac{1}{r^2}\big[\sin\phi\cos\phi-\sin(\phi+\theta)\cos(\phi+\theta)+a^2\sin\theta\cos\theta\big]\\
=&\frac{1}{r^2}\big[a^2\sin\theta\cos\theta-\cos(2\phi+\theta)\sin\theta\big]\\
\geq&\frac{1}{r^2}\big[(1+\ve)\sin\theta-\cos(2\phi+\theta)\sin\theta\big]\\
\geq&\frac{\ve\sin\theta}{r^2}=\frac{\ve}{r^2}\frac{2\mu e^{-at} r^a}{\mu^2+e^{-2at}r^{2a}}\\
\geq&\frac{2\mu e^{-at}\ve}{\mu^2+1}r^{\ve-1}
\end{split}
\eeq
where we have used \eqref{axis6.5}  and \eqref{axis6.3}.
From \eqref{axis6.1} and \eqref{axis6.3}, we  have
\beq\label{axis6.9}
f_t=-\frac{2re^t(\beta+\beta_t)}{e^{2t}\beta^2+r^2}-\frac{2a\mu r^{a}e^{at}}{\mu^2e^{2at}+r^{2a}}
\eeq
and
\beq\label{axis6.10}
f_r=\frac{2\beta e^t}{e^{2t}\beta^2+r^2}+\frac{2a\mu r^{a-1}e^{at}}{\mu^2e^{2at}+r^{2a}}
\eeq
Combining (\ref{axis6.9}) with (\ref{axis6.10}) and using \eqref{axis6.7}, we obtain 
\beq\label{axis6.11}
f_t+rf_r=-\frac{2re^t\beta_t}{e^{2t}\beta^2+r^2}=\frac{2\delta re^{-t}\beta^{\ve}}{e^{2t}\beta^2+r^2}.
\eeq
To prove Claim 3, it suffices to verify the following inequality:
\beq\label{axis6.12}
\frac{2\delta re^{-t}\beta^{\ve}}{e^{2t}\beta^2+r^2}\leq \frac{2\gamma\mu e^{-at}\ve}{\mu^2+1}r^{\ve-1}.
\eeq
Let $s=\frac{r}{e^t\beta}$. Then \eqref{axis6.12} is equivalent to
\beq\label{axis6.13}
\frac{ s^{2-\ve}}{1+s^2}\leq \frac{\gamma\mu \ve}{\delta(\mu^2+1)},\ \forall s>0.
\eeq
It is easy to check that the function $\displaystyle\frac{ s^{2-\ve}}{1+s^2}$ has a maximum $M(\ve)$ depending only on $\ve$. Therefore, if we choose
$$\delta\leq \frac{\gamma\mu \ve}{M(\ve)(\mu^2+1)},$$
then \eqref{axis6.12} holds and hence the Claim 3 follows.

\medskip
\noindent{\bf Claim 4}. {\it For sufficiently large $\mu>0$, there exists $\varphi_0\in C^\infty([0,1])$,
with $\va_0(0)=0$ and $\va_0(1)>\pi$, such that
\beq\label{initial_comp}
f(r,t)\le \varphi_0(r) \ {\rm{on}}\ ([0,1]\times\{0\})\cup (\{0,1\}\times (0,T_0)).
\eeq
}Since $\va_0(1)>\pi$, we can choose a sufficiently large $\mu$ such that $\theta(1,\mu,t)\leq \va_0(1)-\pi$ for any $t\ge 0$.
This, combined with $0<\phi(1,\beta,t)<\pi$, implies
$$
f(1,t)=\phi(1,\beta,t)+\theta(1,\mu,t)\leq \va_0(1), \ \forall 0\le t<T_0.
$$
It is clear that we can find some initial data $\va_0\in C^\infty([0,1])$,
with $\phi_0(0)=0$ and $|\phi_0(1)|>\pi$, such that $\va_0(r)\geq f(r,0)$ for any $r\in[0,1]$.
Hence (\ref{initial_comp}) holds.

It follows from Claim 3, Claim 4, and lemma \ref{comp} that 
$$f\le\va,\  \ {\rm{in}} \ [0,1]\times [0,T_0).$$
This, combined with $f(0,t)=\va(0,t)=0$,  implies that 
$\va_r(0,t)\ge f_r(0,t)$ for $0\le t<T_0$.
Hence, by Claim 2, we have
$$
\va_r(0,t)\rightarrow +\infty,\ \ \mbox{as }\ t\rightarrow T_0^-.
$$
This completes the proof of Theorem \ref{th3}.
\endpf

\begin{remark}
{\rm It is clear that the Claim 3 doesn't hold if the time relexation constant $\gamma=0$.
In fact, it is a very challenging question how to construct solutions with finite time singularity
to (\ref{liquidcrystal2}) when $\gamma=0$. Such an example would shed lights to the question whether 
the equation of viscoelastic flows admits finite time singularity in dimension three.
}
\end{remark}

We finish this section with two remarks on Theorem \ref{blowup}.
\begin{remark}{\rm Under the same assumptions as Theorem \ref{blowup},  does the stress tensor
$\displaystyle\mathcal S:=\nabla d\odot\nabla d-\frac12|\nabla d|^2 \mathbb I_3$ blow up as $t$ approaches
$T_0$?}
\end{remark}

By calculating the stress tensor for the corresponding subsolution, we conjecture that $\mathcal S$
doesn't blow up. Here we sketch the calculation.
Assume $(u,d, P)$ is given by Theorem \ref{blowup}, and $f$ is given by Section 5 above.
Set
$$
\widetilde{d}=\big(\sin f\cos\theta, \sin f\sin \theta, \cos f\big).
$$
Then direct calculations imply
\beq\label{axis7.1}
\nabla \widetilde d\odot\nabla \widetilde d=
\left(
\begin{array}{ccc}
\displaystyle f_r^2\cos\theta^2+\frac{\sin^2f}{r^2}\sin^2\theta
&\displaystyle\big(f_r^2-\frac{\sin^2f}{r^2}\big)\cos\theta\sin\theta\\
\\
\displaystyle \big(f_r^2-\frac{\sin^2f}{r^2}\big)\cos\theta\sin\theta
&\displaystyle f_r^2\sin\theta^2+\frac{\sin^2f}{r^2}\cos^2\theta
\end{array}
\right)
\eeq
and
\beq\label{axis7.2}
|\nabla \widetilde d|^2=f_r^2+\frac{\sin^2f}{r^2}.
\eeq
At $r=0$:  since $\displaystyle\frac{\sin f}{r}=f_r=\frac{2}{\beta(t)}e^{-t}$, it follows that
\beq\label{axis7.3}
\nabla \widetilde d\odot\nabla \widetilde d\ \Big|_{r=0}=
\left(
\begin{array}{ccc}
\frac{4}{\beta^2(t)}
&0\\
\\
0
&\frac{4}{\beta^2(t)}
\end{array}
\right)e^{-2t}
\eeq
and
\beq\label{axis7.4}
|\nabla\widetilde d|^2\Big|_{r=0}=\frac{8}{\beta^2(t)}e^{-2t}.
\eeq
Therefore the stress tensor for $\widetilde d$ at $r=0$ is
\beq\label{axis7.5}
\widetilde{\mathcal S}:=\big(\nabla \widetilde d\odot\nabla \widetilde d
-\frac12|\nabla \widetilde d|^2\mathbb I_{2}\big)\Big|_{r=0}=
\left(
\begin{array}{ccc}
0 &0\\
\\
0&0
\end{array}
\right),
\eeq
where $\mathbb I_2$ is the identity matrix of order $2$.

\begin{remark}{\rm Under the same assumptions as Theorem \ref{blowup}, does the solution
$(u,d,P)$ satisfy the energy dissipation inequality
\beq\label{energy_ineq}
\frac{d}{dt}\Big(\frac12\int_\Omega |u|^2+|\nabla d|^2\Big)
+\int_\Omega \Big(|\nabla u|^2+|\Delta d+|\nabla d|^2 d|^2\Big)\le 0,
\eeq
for $0\le t<T_0$?}
\end{remark}

Since $u(x,y, z, t)=u_0(x,y, z)=(x,y, -2z)$ in $\Omega=B_1^2\times [0,1]$, direct calculations imply
$$\frac{d}{dt}\int_\Omega |u|^2=0, \ \int_\Omega \Delta u\cdot u=0,$$
$$\int_\Omega \langle u\cdot \nabla u, u\rangle=\int_\Omega  u\cdot \nabla \big(\frac{|u|^2}2\big)
=\int_\Omega (x^2+y^2-8z^2)=-\frac{13}{6}\pi,$$
and
$$
\int_\Omega |\nabla u|^2=6\pi.
$$
(\ref{liquidcrystal})$_1$ can be rewritten as
\begin{equation}\label{lcf1}
\nabla\Big(\frac{x^2}2+\frac{y^2}2+2z^2+P\Big)=-\langle\Delta d,\nabla d\rangle.
\end{equation}
Multiplying (\ref{lcf1}) by $u$ and integrating over $\Omega$ yields
\begin{equation}
\frac{13}6 \pi -\int_{\partial\Omega} P u\cdot\nu
=\int_\Omega \langle \Delta d,u\cdot\nabla d\rangle.
\end{equation}
Since the boundary condition of $d$ on $\partial\Omega$ is given by
$$d=d_0 \ {\rm{on}}\ \partial B_1^2\times [0,1], \ \ \frac{\partial d}{\partial\nu}=0 \  {\rm{on}}\ B_1^2\times \{0,1\},$$
multiplying (\ref{liquidcrystal})$_3$  by $\Delta d$ and integrating over $\Omega$ we would have
\begin{equation}
-\frac{d}{dt}\big(\frac12\int_\Omega |\nabla d|^2\big)+\int_\Omega \langle u\cdot\nabla d, \Delta d\rangle
=\int_\Omega\big| |\Delta d+|\nabla d|^2 d\big|^2. 
\end{equation}
Adding (5.24) and (5.25) together, we have
\begin{eqnarray}
&&\frac{d}{dt}\Big(\frac12\int_\Omega |u|^2+|\nabla d|^2\Big)
+\int_\Omega \Big(|\nabla u|^2+|\Delta d+|\nabla d|^2 d|^2\Big)\nonumber\\
&&=6\pi+\frac{13}6 \pi-\int_{\partial\Omega} Pu\cdot {\nu}.
\end{eqnarray}
Therefore (\ref{energy_ineq}) holds iff
\beq\label{bdry_p}
6\pi+\frac{13}6 \pi-\int_{\partial\Omega} Pu\cdot {\nu}\le 0.
\eeq
However, it is unclear whether (\ref{bdry_p}) holds for the solution by Theorem \ref{blowup}.
 
\section{Finite time singularity for generic initial data and proof of Theorem \ref{blow-up2}}
 \setcounter{equation}{0}
\setcounter{theorem}{0}
 
This section is devoted to another construction of finite time singularity of (\ref{liquidcrystal}) for more generic
initial data, in which the solution satisfies the energy dissipation inequality (\ref{energy_ineq}). 
 
First we recall the following result on the existence of local smooth solution to (\ref{liquidcrystal}) and (\ref{IBC}),
whose proof can be found in \cite{lin-lin-wang}.
\begin{lemma}\label{local-solution} For $(u_0, d_0)\in C^\infty_{0,{\rm{div}}}(B_1^3, \mathbb R^3)
 \times C_{\bf e}^\infty(B_1^3, \mathbb S^2)$, there exist $T_0 = T_0(u_0,d_0)> 0$
 and a unique smooth solution $(u,d) \in C^\infty\big(\overline{B_1^3} \times [0,T_0), \mathbb R^3 \times
 \mathbb S^2\big)$ to  the system  (\ref{liquidcrystal}) 
 along with the initial-boundary condition (\ref{IBC}). Moreover, the energy dissipation inequality
 (\ref{energy_ineq}) holds for $0\le t<T_0$.
\end{lemma} 
 
Now we would like to proceed with the proof of Theorem \ref{blow-up2} as follows. 

\smallskip 
\pf   Assume $T_0 > 0$ is the maximal time interval for the short time smooth solution $(u,d)$
by lemma \ref{local-solution}. We want to show

\smallskip
\noindent{\bf Claim 5.} {\it If $\epsilon_0 >0$ is sufficiently small, then $T_0 <1$.}

\smallskip
We argue by contradiction. Suppose that Claim 5 were false. Then for any $\epsilon>0$
we can find $(u_0, d_0) \in C^\infty_{0,{\rm{div}}} (B_1^3, \mathbb R^3) 
\times C^\infty_{\bf e}(B_1^3, \mathbb S^2)$
such that \\
(a) $d_0$ is not homotopic to ${\bf e}$ relative to $\partial B_1^3$, \\
(b) $E(u_0 , d_0) \le \epsilon^2$,\\
and a smooth solution $(u, d) \in C^\infty\big(\overline{B_1^3} \times  [0, 1], \mathbb R^3 \times\mathbb S^2\big)$
 to (\ref{liquidcrystal}) and (\ref{IBC}). Integrating (\ref{energy_ineq}) over $t$ yields that  
 $(u, d)$ satisfies the energy inequality: 
 \begin{equation}\label{energy_ineq1}
 E(u(t), d(t))+\int_0^t\int_{B_1^3} \big(|\nabla u|^2 + |\Delta d +|\nabla d|^2 d|^2\big)\le E(u_0, d_0) 
 \le\epsilon^2,
 \end{equation}
 for all $0\le t\le 1$.
Applying Fubini's theorem to (\ref{energy_ineq1}), we find that there exists $t_1 \in (\frac12, 1)$
such that
\begin{equation}\label{energy_ineq2}
 E(u(t_1), d(t_1))+\int_{B_1^3}\big(|\nabla u(t_1)|^2 + |\Delta d(t_1) + |\nabla d(t_1)|^2d(t_1)|^2\big)\le 8\epsilon^2.
 \end{equation}
From (\ref{energy_ineq2}) and $\epsilon$-apriori estimate (\ref{holder2}) of 
 Theorem \ref{holder1} below, we conclude that there exists a
universal $C > 0$ such that
 \begin{equation}\label{holder3}
 \big[d(t_1)\big]_{C^{\frac12}(B_1^3)} \le C\sqrt{\epsilon}.
 \end{equation} 
Thus $d(t_1)(B_1^3) \subset B^3_{C\sqrt{\epsilon}} ({\bf e}) \cap \mathbb S^2$
 and hence $d(t_1)$ is homotopic to $\bf e$ relative to $\partial B_1^3$, provided 
 $\epsilon >0$ is chosen to be sufficiently small. Since $d\in C^\infty\big(\overline{B_1^3}\times [0, t_1],\mathbb S^2\big)$
 and $d={\bf e}$ on $\partial B_1^3\times [0, t_1]$, we see that $d(t_1)$ is homotopic to $d_0$ relative to 
 $\partial B_1^3$ and hence $d_0$ is homotopic to $\bf e$ relative to 
 $\partial B_1^3$. This contradicts  the assumption (a). 
Hence Claim 5 is true. This completes the proof of Theorem \ref{blow-up2}. \qed\\

Now we need an $\epsilon$-apriori estimate on approximate harmonic maps
from $B_1^3$ to $\mathbb S^2$, which can be proved by suitable modifications of the arguments by Ding-Wang \cite{DW}
 and Lin-Wang \cite{lin-wang}.
\begin{theorem}\label{holder1}
There exist $\epsilon> 0$ and $C>0$ such that if 
$d\in C^\infty_{\bf e}(B_1^3, \mathbb S^2)$ satisfies
\begin{equation}\label{small_energy3}
E(d):=\frac12\int_{B_1^3} |\nabla d|^2 \le\epsilon^2,
\end{equation}
and
\begin{equation}\label{tension_bound}
\int_{B_1^3}|\Delta d+|\nabla d|^2 d|^2\le \epsilon^2.
\end{equation}
Then 
\begin{equation}\label{holder2}
\big[d\big]_{C^\frac12(B_1^3)}\le C\sqrt{\epsilon}.
\end{equation}
\end{theorem}
\pf The proof is based on suitable modifications of that by Ding-Wang \cite{DW} and Lin-Wang \cite{lin-wang}.
For the completeness, we sketch it here. We divide the proof of estimate (\ref{holder2}) into two lemmas.

\begin{lemma}\label{interior_est} Under the same assumptions as Theorem \ref{holder1},
for any fixed $\delta_0\in (\epsilon, 1)$ we have that
\beq\label{interior_holder}
\big[d\big]_{C^\frac12(B_{\frac{r_0}2}(x_0))}\le C\sqrt{\epsilon},
\eeq
holds for any $x_0\in B_{1-\delta_0}^3$ and $0<r_0<\delta_0$.
\end{lemma}

\noindent{\bf Proof of Lemma \ref{interior_est}}: 
In order to show (\ref{interior_holder}), set the tension field of $d$ by
$$\tau(d)\equiv \Delta d+|\nabla d|^2 d: B_1^3\to\mathbb R^3.$$ 
We will first establish a modified energy monotonicity inequality
for approximate harmonic maps in dimension three.
From (\ref{tension_bound}), we see that tension field $\tau(d)\in L^2(B_1^3)$,  and
\begin{equation}\label{tension_bound1}
\big\|\tau(d)\big\|_{L^2(B_1^3)}\le \epsilon.
\end{equation}
Now we have

\smallskip
\noindent{\bf Claim 6}. {\it For any $x_0\in B_1^3$ and $0<r\le R<1-|x_0|$,
it holds
\begin{equation}\label{energy_mono2}
r^{-1}\int_{B_r(x_0)}|\nabla d|^2\le 8 R^{-1}\int_{B_R(x_0)}|\nabla d|^2+8 R\int_{B_R(x_0)}|\tau(d)|^2.
\end{equation}}To see (\ref{energy_mono2}), we assume for simplicity that $x_0=0$. Multiplying the approximate harmonic
map equation
\begin{equation}\label{approx_hm}
\Delta d+|\nabla d|^2 d=\tau(d),   \  \ {\rm{in}}\ \ B_1^3,
\end{equation}
by $x\cdot\nabla d$, integrating the resulting equation over $B_r$ for $0<r<1$, and applying the same argument as
\cite{lin-wang} lemma 5.3, we obtain
\begin{eqnarray}\label{stationary1}
\frac{d}{dr}\Big(r^{-1}\int_{B_r(0)}\big(\frac12|\nabla d|^2-\langle\tau(d), x\cdot\nabla d\rangle\big)\Big)
=r^{-1}\int_{\partial B_r(0)}\big|\frac{\partial d}{\partial r}\big|^2-\int_{\partial B_r(0)}\big\langle \tau(d),
\frac{\partial d}{\partial r}\big\rangle.
\end{eqnarray}
Thus (\ref{energy_mono2}) follows by integrating (\ref{stationary1}) over $r$ and H\"older's inequality.

It follows from (\ref{small_energy3}), (\ref{tension_bound}), and (\ref{energy_mono2}) that 
\begin{eqnarray}
\label{small_renormal}
&&\sup\Big\{r^{-1}\int_{B_r(x)}|\nabla d|^2: \ x\in B_{1-\delta_0}^3, \ 0<r\le \delta_0\Big\}\nonumber\\
&&\le 8\delta_0^{-1} \int_{B_1^3}|\nabla d|^2+ 8\delta_0\int_{B_1^3}|\tau(d)|^2\le C\epsilon^2(1+\delta_0^{-1})
\le C\epsilon.
\end{eqnarray}
Next we have

\smallskip
\noindent{\bf Claim 7}. {\it For any $\alpha\in (0, \frac12]$,  there exists $\theta_0\in (0,\frac12)$ such that  for
any $x_0\in B_{1-\delta_0}^3$ and $0<r<\delta_0$ it holds
\beq\label{renormal_decay1}
(\theta_0 r)^{-1}\int_{B_{\theta_0 r}(x_0)}|\nabla d|^2\le \theta_0^{2\alpha} r^{-1}\int_{B_r(x_0)}|\nabla d|^2
+C\theta_0^{-1} r \int_{B_r(x_0)}|\tau(d)|^2.
\eeq}
To show (\ref{renormal_decay1}), first recall the Morrey space $M^{2,2}(U)$ for $U\subset\R^3$:
$$M^{2,2}(U):=\Big\{f\in L^2_{\rm{loc}}(U)\ | \ \|f\|_{M^{2,2}(U)}^2\equiv
\sup_{B_r(x)\subset U} r^{-1}\int_{B_r(x)}|f|^2<+\infty\Big\}.$$
From (\ref{small_renormal}),  we have that $\nabla d\in M^{2,2}(B_{\delta_0}(x))$ for any $x\in B_{1-\delta_0}^3$,
and
\beq \label{small_renormal1}
\big\|\nabla d\big\|_{M^{2,2}(B_{\delta_0}(x))}\le C\epsilon, \ \forall x\in B_{1-\delta_0}^3.
\eeq
For $x_0\in B_{1-\delta_0}^3$ and $0<r\le 1-\delta_0$, 
let $\widetilde{d}\in H^1(\mathbb R^3,\mathbb R^3)$ be an extension of $d$  such that
$\widetilde{d}=d$ in $B_r(x_0)$,  $|\widetilde{d}(x)|\le 2$ for $x\in\R^3$, and
\beq\label{extension}\begin{cases}
\big\|\nabla\widetilde d\big\|_{L^2(\R^3)}\leq C \big\|\nabla d\big\|_{L^2(B_r(x_0))},\\
\big\|\nabla \widetilde{d}\big\|_{M^{2,2}(\R^3)}\leq C \big\|\nabla d\big\|_{M^{2,2}(B_r(x_0))}.
\end{cases}
\eeq
By the Helmholtz decomposition, there exist $G\in H^1(\mathbb R^3)$ and $H\in L^2(\mathbb R^3, \R^3)$ such that
\beq\label{H-decom}
\begin{cases}
\nabla\widetilde d\times \widetilde d=\nabla G+H, \ \nabla\cdot H=0\ \ {\rm{in}}\ \R^3,\\
\big\|\nabla G\big\|_{L^2(\R^3)}+\big\|H\big\|_{L^2(\R^3)}\leq C \big\|\nabla\widetilde d\big\|_{L^2(\R^3)}
\leq C \big\|\nabla d\big\|_{L^2(B_r(x_0))}.
\end{cases}
\eeq
Applying the Poincar\'e inequality, the duality between the Hardy space $\mathcal H^1(\R^3)$ and BMO($\R^3$) (see \cite{evans}, \cite{helein},
or \cite{lin-wang}), we can estimate 
\begin{eqnarray}\label{H-est}
\int_{\R^3}|H|^2&=&\int_{\R^3}H\cdot(\nabla\widetilde d\times\widetilde d-\nabla G)
=\int_{\R^3} H\cdot\nabla\widetilde d\times \widetilde d\nonumber\\
&\leq& C\big\|H\cdot\nabla \widetilde d\big\|_{\mathcal H^1(\R^3)}\big\|\widetilde d\big\|_{{\rm{BMO}}(\R^3)}\nonumber\\
&\leq& C\big\|H\big\|_{L^2(\R^3)}\big\|\nabla\widetilde d\big\|_{L^2(\R^3)}\big\|\nabla \widetilde d\big\|_{M^{2,2}(\R^3)}\nonumber\\
&\leq& C\big\|\nabla d\big\|_{L^2(B_r(x_0))}^2\big\|\nabla d\big\|_{M^{2,2}(B_r(x_0))}
\leq C \epsilon \big\|\nabla d\big\|_{L^2(B_r(x_0))}^2.
\end{eqnarray}
From (\ref{approx_hm}),  $G$ solves
\beq\label{approx-hm1}
\Delta G=\tau(d)\times d \ \ {\rm{in}}\ \ B_r(x_0).
\eeq
By the standard $L^2$-estimate, we have that for any $0<\theta<1$,
\beq\label{G-est}
\int_{B_{\theta r}(x_0)}|\nabla G|^2\leq C \theta^3 \int_{B_r(x_0)}|\nabla d|^2+ r^2\int_{B_r(x_0)}|\tau(d)|^2.
\eeq
Putting (\ref{H-est}) and (\ref{G-est}) together, we obtain
\beq\label{grad-est}
\frac{1}{\theta r}\int_{B_{\theta r}(x_0)}|\nabla d|^2
\le C(\theta^2+\theta^{-1}\epsilon) \frac{1}{r}\int_{B_r(x_0)}|\nabla d|^2+C\theta^{-1} r\int_{B_r(x_0)}|\tau(d)|^2,
\eeq
for any $x_0\in B_{1-\delta_0}^3$, $0<r<\delta_0$, and $0<\theta<1$.

For any $0<\alpha\le \frac12$, first choosing $\theta=\theta_0\in (0,1)$ such that $C\theta_0^2\le \frac12\theta_0^{2\alpha}$
and then choosing $\epsilon=\epsilon_0>0$ such that $C\theta_0^{-1}\epsilon_0\le \frac12\theta_0^{2\alpha}$, we obtain
(\ref{renormal_decay1}). It is standard that iterations of (\ref{renormal_decay1}) imply that
\beq \label{morrey_decay1}
\frac{1}{s}\int_{B_s(x_0)}|\nabla d|^2\le \big(\frac{s}{r})^{2\alpha}\int_{B_r(x_0)}|\nabla d|^2+Cs\int_{B_r(x_0)}|\tau(d)|^2,
\eeq
holds for any $x_0\in B_{1-\delta_0}^3$ and $0<s\le r<\delta_0.$ It is clear that (\ref{interior_holder})
follows from (\ref{morrey_decay1}) and Morrey's
decay lemma (cf. \cite{morrey}). \qed

\begin{lemma} \label{bdry_est} Under the same assumptions as in Theorem \ref{holder1}, 
there exists $0<\delta_0\le \frac12$ depending only on
$\partial B_1^3$ such that for any $x_0\in\partial B_1^3$ and $0<r_0\le 2\delta_0$, it holds 
\beq\label{bdry_holder}
\big[d\big]_{C^\frac12(B_{\frac{r_0}2}^+(x_0))}\le C\sqrt{\epsilon},
\eeq
where $B_r^+(x_0):=B_r(x_0)\cap B_1^3$ for $r>0$.
\end{lemma}

\noindent{\bf Proof of Lemma \ref{bdry_est}}: The strategy to show (\ref{bdry_holder}) is similar to
lemma \ref{interior_est}, that is to establish a modified boundary monotonicity inequality and a boundary
energy decay property. More precisely, we first need

\bigskip
\noindent{\bf Claim 8}. {\it There exists  $0<\delta_0\le \frac12$, depending only on $\partial B_1^3$,
such that for any $x_0\in \partial B_1^3$ and $0<r\le R\le \delta_0$, it holds
\begin{equation}\label{energy_mono3}
r^{-1}\int_{B_r^+(x_0)}|\nabla d|^2\le 8 R^{-1}\int_{B_R^+(x_0)}|\nabla d|^2+16 R\int_{B_R^+(x_0)}|\tau|^2.
\end{equation}}

For $r>0$, set $\mathbb B_r^+=B_r(0)\cap\big\{x=(x_1,x_2, x_3)\in\R^3:\ x_3\ge 0\big\}$, 
$T_r=\mathbb B_r^+\cap\big\{x=(x_1,x_2,x_3)\in\R^3: \ x_3=0\big\}$, and $S_r^+=\partial\mathbb B_r^+\setminus T_r$.
By the standard boundary flatten argument, there exists $\delta_0>0$ depending only on $\partial B_1^3$
such that for any $x_0\in\partial B_1^3$ there exists a smooth diffeomorphism $\Phi_0: \mathbb B_{\delta_0}^+
\to B_{\delta_0}^+(x_0)$ such that $\Phi_0(0)=x_0$, $\Phi_0(T_{\delta_0})=B_{\delta_0}^+(x_0)\cap\partial B_1^3$, and
\beq\label{bdry_flat}
|\Phi_0(x)-x_0|\le C|x|^2, \ |\nabla\Phi_0(x)-{\rm{Id}}|\le C|x|, |\nabla^2\Phi_0(x)|\le C,
\ \forall x\in \mathbb B_{\delta_0}^+,
\eeq
where Id denotes the identity map in $\R^3$.  

Set $g=\Phi_0^*(g_0)$ in $\mathbb B_{\delta_0}^+$, pull-back of the euclidean metric $g_0=dx^2$ in $B_{\delta_0}^+(x_0)$.
Consider $\widehat{d}=d\circ\Phi_0: \mathbb B_{\delta_0}^+\to\mathbb S^2$. Then it is not hard to check that
$\widehat{d}: (\mathbb B_{\delta_0}^+, g)\to \mathbb S^2$ is an approximated harmonic map, i.e.,
\beq\label{approx_hm2}
\Delta_g \widehat d+\big|\nabla \widehat d\big|_g^2\widehat{d}=\tau\big(\widehat{d}\ \big) \ \ {\rm{in}}\ \ \mathbb B_r^+;
\ \widehat{d}={\bf e} \ \ {\rm{on}}\ \ T_{\delta_0},
\eeq
and
\beq\label{tension_bd}
\big|\tau\big(\widehat{d}\ \big)\big|(x)\le C|\tau(d)|(\Phi_0(x)), \ \forall x\in \mathbb B_{\delta_0}^+,
\eeq
where $\Delta_g$ is the Laplace operator with respect to $g$, and $|\nabla\widehat{d}|_g^2$ is the Dirichlet
energy density of $\widehat d$ with respect to $g$. Note that by (\ref{bdry_flat}), $g$ satisfies
\beq\label{approx_flat}
|g(x)-g_0|\le C|x|^2, \ |\nabla g(x)|\le C|x|, \ |\nabla^2 g(x)|\le C, \ \forall x\in \mathbb B_{\delta_0}^+.
\eeq
Based on (\ref{approx_flat}), (\ref{approx_hm2}), and (\ref{tension_bd}), we may for simplicity assume 
that $x_0=0$ and $d: (\mathbb B_{\delta_0}^+, g_0)\to \mathbb S^2$ 
is an approximated harmonic map with tension field $\tau(d)$. Since $d={\bf e}$ on $T_{\delta_0}$,
it is easy to see $x\cdot\nabla d=0$ on $T_{\delta_0}$. For $0<r\le\delta_0$,
multiplying (\ref{approx_hm}) by $x\cdot\nabla d$ and integrating the resulting equation over $\mathbb B_r^+$, we get
\begin{eqnarray*}
\frac12\int_{\mathbb B_r^+}|\nabla d|^2-\frac{r}{2}\int_{S_r^+}|\nabla d|^2
+r\int_{S_r^+}\big|\frac{\partial d}{\partial r}\big|^2=\int_{\mathbb B_r^+}\langle\tau(d), x\cdot\nabla d\rangle.
\end{eqnarray*}
This implies
\begin{eqnarray}\label{bdry_mono2}
\frac{d}{dr}\Big(r^{-1}\int_{\mathbb B_r^+}\big(\frac12|\nabla d|^2-\langle\tau(d), x\cdot\nabla d\rangle\big)\Big)
=r^{-1}\int_{S_r^+}\big|\frac{\partial d}{\partial r}\big|^2-\int_{S_r^+}\big\langle \tau(d),
\frac{\partial d}{\partial r}\big\rangle.
\end{eqnarray}
Integrating (\ref{bdry_mono2}) over $0<r<R\le\delta_0$ and applying H\"older's inequality, we obtain
$$
r^{-1}\int_{\mathbb B_r^+}|\nabla d|^2\le 8R^{-1}\int_{\mathbb B_R^+}|\nabla d|^2+16 R\int_{\mathbb B_R^+}|\tau(d)|^2.
$$
This gives (\ref{energy_mono3}).

Next we need

\noindent {\bf Claim 9}. {\it For any $\alpha\in (0, \frac12)$,  there exists $\theta_0\in (0,\frac12)$ such that  for
any $x_0\in \partial B^3$ and $0<r<\delta_0$, it holds
\beq\label{renormal_decay2}
(\theta_0 r)^{-1}\int_{B_{\theta_0 r}^+(x_0)}|\nabla d|^2\le C\theta_0^{2\alpha} r^{-1}\int_{B_r^+(x_0)}|\nabla d|^2
+C\theta_0 r \int_{B_r^+(x_0)}|\tau(d)|^2.
\eeq}
For simplicity, we again assume $x_0=0$ and $B_{\delta_0}^+(x_0)=\mathbb B_{\delta_0}^+$.
The proof is similar to \cite{DW} lemma 3.3, and we only sketch it.
To obtain (\ref{renormal_decay2}), we perform
suitable extensions of $d$ to $\mathbb B_{\delta_0}$ as follows.
Let $\widetilde{d}: \mathbb B_{\delta_0}$ be the extension of $d$ that is even
with respect to $x_3$, and let $\widetilde {\tau}: \mathbb B_{\delta_0}\to \R^3$ be the extension
of $\tau(d)$ that is odd with respect to $x_3$. Define $w: \mathbb B_{\delta_0}\to \R^3$ by
\beq\label{odd_ext}
w(x)=\begin{cases} (d-{\bf e})(x) & {\rm{if}}\  x=(x_1,x_2,x_3)\in \mathbb B_{\delta_0}
\ {\rm{and}}\ x_3\ge 0,\\
-(d-{\bf e})(x_1, x_2, -x_3) & {\rm{if}}\ x=(x_1,x_2,x_3)\in \mathbb B_{\delta_0}
\ {\rm{and}}\ x_3<0.
\end{cases}
\eeq

For $0<r\le\delta_0$, let $\eta\in C_0^\infty(\mathbb R^3)$ be even with respect to
$x_3$, $0\le \eta\le 1$, $\eta=1$ in $\mathbb B_{\frac{r}2}$, $\eta=0$ outside $\mathbb B_r$,
and $|\nabla\eta|\le Cr^{-1}.$  To proceed with the proof, we need 
\beq\label{div1}
\nabla\cdot\big(\nabla w\times \widetilde{d}\ \big)=\widetilde{\tau}\times \widetilde{d} 
\ \ \ {\rm{in}}\ \ \mathcal D'(\mathbb B_{\delta_0}).
\eeq
To see (\ref{div1}), let $\phi\in C_0^\infty(\mathbb B_{\delta_0})$ and write $\phi=\phi_{\rm{e}}+\phi_{\rm{o}}$,
here $\phi_{\rm{e}}$ ($\phi_{\rm{o}}$) is even (odd) with respect to $x_3$ respectively. 
Then (\ref{div1}) follows from
\begin{eqnarray*}
\int_{\mathbb B_{\delta_0}} \nabla w\times\widetilde{d}\cdot\nabla \phi
&=&2\int_{\mathbb B_{\delta_0}^+} \nabla w\times\widetilde{d}\cdot\nabla \phi_{\rm{o}}
=-2\int_{\mathbb B_{\delta_0}^+} \nabla\cdot(\nabla d\times d)\cdot \phi_{\rm{o}}\\
&=&-2\int_{\mathbb B_{\delta_0}^+} \tau(d)\times d\cdot \phi_{\rm{o}}
=-\int_{\mathbb B_{\delta_0}} \widetilde{\tau}\times \widetilde{d}\cdot \phi.
\end{eqnarray*}
It follows from (\ref{energy_mono2}) and (\ref{energy_mono3}) that  $\nabla w\in M^{2,2}(\mathbb B_{\delta_0})$, and
\beq\label{morrey01}
\big\|\nabla w\big\|_{M^{2,2}(\mathbb B_{\delta_0})}\le C\epsilon.
\eeq
From (\ref{div1}), it is easy to check that 
$\nabla\cdot\big(\eta^2\nabla w\times\widetilde{d}\ \big)\in L^2(\R^3)$, and
\beq\label{l2-div}
\Big\|\nabla\cdot\big(\eta^2\nabla w\times\widetilde{d}\ \big)\Big\|_{L^2(\R^3)}
\le C\Big(r^{-1}\|\nabla d\|_{L^2(\mathbb B_r^+)}+\|\tau(d)\|_{L^2(\mathbb B_r^+)}\Big).
\eeq
Observe that 
\begin{eqnarray*}
&&2\int_{\mathbb B_{\frac{r}2}^+}|\nabla d|^2\le 2\int_{\mathbb R^3_+}\eta^2|\nabla w\times {d}|^2
=\int_{\R^3}\eta^2\big|\nabla w\times \widetilde{d}\big|^2\\
&&=-\int_{\R^3}\nabla(\eta^2\nabla w\times\widetilde{d}\ )\cdot(w\times\widetilde{d}\ )
+\int_{\R^3}\eta^2[(\nabla w\times \widetilde{d})\times \nabla\widetilde{d}-\lambda]\cdot w
+\lambda\int_{\R^3}\eta^2 w\\
&&=I+II+III,
\end{eqnarray*}
where
$$\lambda=\frac{\int_{\R^3}\eta^2 (\nabla w\times \widetilde{d})\times\nabla \widetilde{d}}{\int_{\R^3}\eta^2}.$$
By the Poincar\'e inequality, we have
$$|III|\le |\lambda|\int_{\mathbb B_r^+}|w|
\le C\big(r^{-1}\int_{\mathbb B_r^+}|\nabla d|^2\big)^{\frac12}\int_{\mathbb B_r^+}|\nabla d|^2.
$$
Applying (\ref{l2-div}),  $I$ can be estimated by
\begin{eqnarray*}
|I|&\le& \Big\|\nabla\cdot\big(\eta^2\nabla w\times\widetilde{d}\ \big)\Big\|_{L^2(\R^3)}\Big\|w\Big\|_{L^2(\mathbb B_r^+)}\\
&\le& C\Big(r^{-1}\|\nabla d\|_{L^2(\mathbb B_r^+)}+\|\tau(d)\|_{L^2(\mathbb B_r^+)}\Big)\Big\|d-{\bf e}\Big\|_{L^2(\mathbb B_r^+)}\\
&\le& \gamma \Big(\|\nabla d\|_{L^2(\mathbb B_r^+)}^2+r^2\|\tau(d)\|_{L^2(\mathbb B_r^+)}^2\Big)
+C\gamma^{-1} r^{-2}\int_{\mathbb B_r^+}|d-{\bf e}|^2
\end{eqnarray*}
for any $\gamma\in (0,1)$.
Let $\widetilde{w}:\mathbb R^3\to\mathbb R^3$ be an extension of $w$ such that 
$$\big\|\nabla\widetilde{w}\big\|_{M^{2,2}(\R^3)}\le C\big\|\nabla w\big\|_{M^{2,2}(B_{\delta_0})}.$$
Then, from (\ref{morrey01}) and the Poincar\'e inequality, we have
\beq\label{morrey02}
\big[\widetilde{w}\big]_{{\rm{BMO}}(\mathbb R^3)}\le C\big\|\nabla \widetilde{w}\big\|_{M^{2,2}(\mathbb R^3)}
\le C\big\|\nabla w\big\|_{M^{2,2}(\mathbb B_{\delta_0})}\le C\sqrt{\epsilon}.
\eeq
We use the duality between $\mathcal H^1(\R^3)$ and BMO($\R^3$), similar to \cite{DW} lemma 3.3 and lemma 2.6, 
to estimate $II$ by
\begin{eqnarray*}
\big|II\big|&\le& C\Big\|\eta^2[(\nabla w\times \widetilde{d})\times \nabla\widetilde{d}-\lambda]\Big\|_{\mathcal H^1(\R^3)}
\big[\widetilde{w}\big]_{{\rm{BMO}}(\mathbb R^3)}\\
&\le& C\Big[\big\|\nabla d\big\|_{L^2(\mathbb B_r^+)}^2+r^2\big\|\nabla\cdot(\nabla d\times d)\big\|_{L^2(\mathbb B_r^+)}^2\Big]
\big\|\nabla w\big\|_{M^{2,2}(\mathbb B_{\delta_0})}\\
&\le& C\sqrt{\epsilon} \Big[\big\|\nabla d\big\|_{L^2(\mathbb B_r^+)}^2+r^2\big\|\tau(d)\big\|_{L^2(\mathbb B_r^+)}^2\Big].
\end{eqnarray*}
Putting these estimates together, we obtain
\beq\label{renormal_decay5}
\frac{1}{r/2}\int_{\mathbb B_{r/2}^+}|\nabla d|^2
\le C(\sqrt{\epsilon}+\gamma) \Big[\frac{1}{r}\int_{\mathbb B_r^+}|\nabla d|^2
+r\int_{\mathbb B_r^+}|\tau(d)|^2\Big]
+C\gamma^{-1}\frac{1}{r^{3}}\int_{\mathbb B_r^+}|d-{\bf e}|^2.
\eeq
To deduce (\ref{renormal_decay2}) from (\ref{renormal_decay5}), we need 

\smallskip
\noindent{\bf Claim 10}. {\it There exists $\theta_0=\theta_0(\epsilon)$ such that for $0<r<\delta_0$,
\beq\label{mean_decay}
\frac{1}{(\theta_0r)^{3}}\int_{\mathbb B_{\theta_0 r}^+}|d-{\bf e}|^2\le C\theta_0^2 \frac{1}{r}\int_{\mathbb B_r^+}|\nabla d|^2.
\eeq}
We prove (\ref{mean_decay}) by contradiction. Suppose that it were false. Then for any $\theta\in (0,1)$, there exist
$\epsilon_k\rightarrow 0$, $0<r_k<\delta_0$,  ${\bf e}_k\in \mathbb S^2$, 
and $d_k\in C^\infty\big(\overline{\mathbb B_{\delta_0}^+},\mathbb S^2\big)$ such that $d_k={\bf e}_k$ on $T_{\delta_0}$,
\beq\label{renormal_small3}
\displaystyle r_k^{-1}\int_{\mathbb B_{r_k}^+}|\nabla d_k|^2=\epsilon_k^2
\ {\rm{and}}\  r_k\int_{\mathbb B_{r_k}^+}|\tau(d_k)|^2\le\epsilon_k^2,
\eeq
but
\beq\label{mean_no_decay}
\frac{1}{(\theta r_k)^{3}}\int_{\mathbb B_{\theta r_k}^+}|d_k-{\bf e}_k|^2>k
\theta^2 \frac{1}{r_k}\int_{\mathbb B_{r_k}^+}|\nabla d_k|^2.
\eeq
Define $\displaystyle\widetilde{d_k}(x)=\frac{d_k(r_k x)-{\bf e}_k}{\epsilon_k}: \mathbb B_{1}^+\to \mathbb S^2$.
Then we have 
\beq\label{renormal_small4}
\widetilde{d_k}\big|_{T_1}=0,\ 
\int_{\mathbb B_1^+}|\nabla\widetilde{d_k}|^2=1,,
\eeq
and
\beq\label{approx_hm3}
\Delta \widetilde{d_k}=f_k:=-\epsilon_k|\nabla\widetilde{d_k}|^2 d_k+\epsilon_k^{-1}\tau(d_k)\rightarrow 0
\ {\rm{in}}\  L^2(\mathbb B_1^+),
\eeq
but
\beq \label{mean_nodecay}
\frac{1}{\theta^3}\int_{\mathbb B_\theta^+}|\widetilde{d_k}|^2>k\theta^2.
\eeq
From (\ref{renormal_small4}), $\big\{\widetilde{d_k}\big\}\subset H^1(\mathbb B_1^+, \mathbb R^3)$ is bounded. We may
assume that $\widetilde{d_k}\rightarrow \widetilde{d}$ weakly in $H^1(\mathbb B_1^+)$, strongly in
$L^2(\mathbb B_1^+)$. It follows from (\ref{renormal_small4}) and (\ref{approx_hm3}) that
\beq\label{renormal_small5}
\widetilde{d}\big|_{T_1}=0,\ 
\int_{\mathbb B_1^+}|\nabla\widetilde{d}|^2\le1,
\eeq
and
\beq\label{approx_hm4}
\Delta \widetilde{d}=0 \ \ {\rm{in}}\ \ \mathbb B_1^+.
\eeq
By the standard theory of harmonic functions, we have that for any $\theta\in (0,1)$,
\beq\label{decay_harmonic}
\frac{1}{\theta^3}\int_{\mathbb B_\theta^+}|\widetilde{d}|^2\le C\theta^2.
\eeq
Since $\widetilde{d_k}\rightarrow \widetilde{d}$ in $L^2(\mathbb B_1^+)$, (\ref{decay_harmonic})
contradicts to (\ref{mean_nodecay}). Hence Claim 10  is proven.

Putting (\ref{mean_decay}) into (\ref{renormal_decay5}), we can obtain that 
\begin{eqnarray}\label{renormal_decay6}
\frac{1}{\theta_0r}\int_{\mathbb B_{\theta_0r}^+}|\nabla d|^2
&\le& C(\sqrt{\epsilon}+\gamma) \Big[\frac{1}{2\theta_0r}\int_{\mathbb B_{2\theta_0r}^+}|\nabla d|^2
+\theta_0r\int_{\mathbb B_{2\theta_0r}^+}|\tau(d)|^2\Big]\nonumber\\
&&+C\gamma^{-1}\theta_0^2\frac{1}{r^{3}}\int_{\mathbb B_r^+}|\nabla d|^2\nonumber\\
&\le& C\big[(\sqrt{\epsilon}+\gamma)\theta_0^{-1}+\gamma^{-1}\theta_0^2\big]\frac{1}{r^{3}}\int_{\mathbb B_r^+}|\nabla d|^2
+C\theta_0r\int_{\mathbb B_{r}^+}|\tau(d)|^2\nonumber\\
&\le& C\theta_0^{2\alpha}\frac{1}{r^{3}}\int_{\mathbb B_r^+}|\nabla d|^2
+C\theta_0r\int_{\mathbb B_{r}^+}|\tau(d)|^2,
\end{eqnarray}
for any $0<\alpha<\frac12$, provided we choose $\gamma=\theta_0^{2-2\alpha}$ and $\epsilon\le\theta_0^{2+4\alpha}$.
Hence (\ref{renormal_decay2}) is proven.

It is standard that iterations of (\ref{renormal_decay2}) imply that for any $\alpha\in (0,\frac12)$, 
\beq\label{morrey_decay2}
\frac{1}{s}\int_{B_s^+(x_0)}|\nabla d|^2\le C\big(\frac{s}{r}\big)^{2\alpha}\int_{B_r^+(x_0)}|\nabla d|^2
+Cs\int_{B_r^+(x_0)}|\tau(d)|^2
\eeq
holds for any $x_0\in\partial B_1^3$ and $0<s<r<\delta_0$. Combining (\ref{morrey_decay2}) with
(\ref{morrey_decay1}), we can obtain that for any $0<\alpha<\frac12$, 
\beq\label{morrey_decay3}
\frac{1}{s}\int_{B_s(x_0)\cap B_1^3}|\nabla d|^2\le C\big(\frac{s}{r}\big)^{2\alpha}\int_{B_r(x_0)\cap B_1^3}|\nabla d|^2
+Cs\int_{B_r(x_0)\cap B_1^3}|\tau(d)|^2
\eeq
holds for any $x_0\in \overline{B_1^3}$ and $0<s<r<\delta_0$.  This, combined with Morrey's decay lemma (see \cite{morrey})
implies that for any $\alpha\in (0,\frac12)$,
\beq\label{global_holder}
\big[d\big]_{C^\alpha\big(\overline{B_1^3}\big)}\le C\Big[\big\|\nabla d\big\|_{L^2(B_1^3)}
+\big\|\tau(d)\big\|_{L^2(B_1^3)}\Big].
\eeq
It is well known that we can use the equation (\ref{approx_hm2}) to
improve the estimate (\ref{global_holder}) to the case that $\alpha=\frac12$. Hence lemma \ref{bdry_est} is proven.
This completes the proof of Theorem 6.2. Thus Theorem 1.3 is proven. 
\qed

\bigskip
\noindent{\bf Acknowledgements}. The first author is partially supported by NSF grants
DMS 1412005 and DMS 1159937. The second author is partially supported by NSF grants DMS 1065964 and
DMS 1159313. The third author is partially supported by NSF grants DMS 1412005, DMS 1216938 and DMS 1159937.
The fourth author is partially supported by NSF grant DMS 1522869 and NSFC grant 11128102.

\end{document}